\documentclass[a4paper, english]{smfart} 

\usepackage{placeins} 

\usepackage{setspace}
\usepackage[utf8]{inputenc}

\usepackage[british]{babel}

\usepackage[OT2,T1]{fontenc}
\usepackage{amssymb} 
\usepackage{mathrsfs} 

\usepackage{stmaryrd}
\usepackage{amsthm}

\usepackage{amsmath}

\allowdisplaybreaks[1]

\usepackage{graphicx}
\usepackage{manfnt} 

 \usepackage{booktabs} 

\usepackage[margin=3cm]{geometry}

\usepackage[all]{xy}

\usepackage{math tools} 

\usepackage{enumitem}

\usepackage[x11names]{xcolor}

\usepackage[pagebackref=true, colorlinks=true, linkcolor=black, citecolor=black, urlcolor=black]{hyperref}
\usepackage{smfthm}

\usepackage[msc-links, alphabetic]{amsrefs}

\DeclareSymbolFont{cyrletters}{OT2}{wncyr}{m}{n}
\DeclareMathSymbol{\Sha}{\mathalpha}{cyrletters}{"58}

\newtheorem{theoA}{Theorem}

\newtheorem*{coro*}{Corollary}
\newtheorem*{conj*}{Conjecture}
\newtheorem*{lemm*}{Lemma}

\providecommand{\twomat}[4]{\left(\begin{array}{cc}#1&#2\\#3&#4\end{array}\right)}
\providecommand{\smalltwomat}[4]{\left(\begin{smallmatrix}#1&#2\\#3&#4\end{smallmatrix}\right)}

\theoremstyle{definition}

\theoremstyle{remark}
\newtheorem{remark*}{Remark}

\NumberTheoremsIn{subsection}

\numberwithin{equation}{subsection}
\numberwithin{table}{subsection}

\newcommand{\ot}{\otimes}

\newcommand{\ts}{\times}
\newcommand{\cd}{\cdot}

\newcommand{\beq}{\begin{equation}\begin{aligned}}
\newcommand{\eeq}{\end{aligned}\end{equation}}

\newcommand{\beqq}{\begin{equation*}\begin{aligned}}
\newcommand{\eeqq}{\end{aligned}\end{equation*}}

\newcommand{\lb}[1]{\label{#1}}

\newcommand{\nek}{Nekov{\'a}{\v{r}}}

\newcommand{\one}{\mathbf{1}}
\newcommand{\Q}{\mathbf{Q}}
\newcommand{\qqq}{\mathbf{q}}

\newcommand{\GL}{\mathrm{GL}}
\newcommand{\G}{\mathrm{G}}
\renewcommand{\H}{\mathrm{H}}

\newcommand{\R}{\mathbf{R}}
\newcommand{\Z}{\mathbf{Z}}

\newcommand{\frakm}{\mathfrak{m}}

\newcommand{\into}{\hookrightarrow}

\newcommand{\lan}{\langle}
\newcommand{\ran}{\rangle}

\newcommand{\tht}{\theta}
\newcommand{\lm}{\lambda}
\newcommand{\Lm}{\Lambda}
\newcommand{\sg}{\sigma}

\newcommand{{\calG}}{\mathscr{G}}

\newcommand{\bC}{\mathbf{C}}

\newcommand{\OO}{\mathscr{O}}
\newcommand{\A}{\mathbf{A}}

\newcommand{\bks}{\backslash}

\newcommand{\eps}{\varepsilon}
\newcommand{\vphi}{\varphi}
\newcommand{\vpi}{\varpi}

\newcommand{\wtil}{\widetilde}
\newcommand{\Res}{\mathrm{Res}}

\newcommand{\ord}{\mathrm{ord}} 

\newcommand{\vol}{\mathrm{vol}}
\newcommand{\Tr}{\mathrm{Tr}}

\newcommand{\Gal}{\mathrm{Gal}}

\newcommand{\Fr}{\mathrm{Fr}}

\newcommand{\Ch}{\mathrm{Ch}}

\newcommand{\Hom}{\mathrm{Hom}\,}
\newcommand{\End}{\mathrm{End}\,}

\newcommand{\Spec}{\mathrm{Spec}\,}

\newcommand{\id}{\mathrm{id}}

\renewcommand{\r}[1]{\mathrm{#1}}
\newcommand{\s}[1]{\mathscr{#1}}
\renewcommand{\sf}[1]{\mathsf{#1}}
\renewcommand{\(}{\left(}
\renewcommand{\)}{\right)}

\newcommand{\ol}[1]{\overline{#1}{}}
\newcommand{\wt}[1]{\widetilde{#1}{}}
\newcommand{\ul}{\underline}
\renewcommand{\geq}{\geqslant}

\newcommand{\rF}{\r F}
\newcommand{\rG}{\r G}
\newcommand{\rH}{\r H}

\newcommand{\rM}{\r M}
\newcommand{\rN}{\r N}

\newcommand{\rP}{\r P}

\newcommand{\rT}{\r T}
\newcommand{\rU}{\r U}

\newcommand{\ra}{\r a}

\newcommand{\rc}{\r c}
\newcommand{\rd}{\,\r d}

\newcommand{\ri}{\r i}

\renewcommand{\rm}{\r m}

\newcommand{\rp}{\r p}

\newcommand{\rr}{\r r}

\newcommand{\rt}{\r t}
\newcommand{\ru}{\r u}

\newcommand{\rw}{\r w}

\newcommand{\sC}{\s C}

\newcommand{\sF}{\s F}

\newcommand{\sH}{\s H}

\newcommand{\sL}{\s L}
\newcommand{\sM}{\s M}

\newcommand{\sO}{\s O}

\newcommand{\sS}{\s S}

\newcommand{\sV}{\s V}
\newcommand{\sW}{\s W}
\newcommand{\sX}{\s X}

\newcommand{\Qpb}{\ol{\Q}_{p}}
\newcommand{\Zpb}{\ol{\Z}_{p}}

\newcommand{\cyc}{\textstyle\circ}

\usepackage{wasysym}

\DeclareMathOperator{\Span}{Span}

\newsavebox\tempbox
\let\svwidetilde\widetilde
\renewcommand\widetilde[1]{\sbox\tempbox{$#1$}\svwidetilde{\usebox{\tempbox}}}

   \def\XXint#1#2#3{{\setbox0=\hbox{$#1{#2#3}{\int}$}
        \vcenter{\hbox{$#2#3$}}\kern-.5\wd0}}

   \newcommand{\Herm}{\r{Herm}}

\title{Euler systems for conjugate-symplectic motives}
\author{Daniel Disegni} 

\address{Department of Mathematics, Ben-Gurion University of the Negev, Be'er Sheva 84105, Israel}
\address{Aix-Marseille University, CNRS, I2M - Institut de Math\'ematiques de Marseille, campus de Luminy, 13288 Marseille, France
}
\email{daniel.disegni@univ-amu.fr}

\thanks{Research partly supported by ISF grant 1963/20 and BSF grant 2018250.}

\begin{document}

\maketitle

\begin{abstract} 
Consider a conjugate-symplectic geometric representation $\rho$ of the Galois group of a CM field. Under the assumption that $\rho$ is automorphic, even-dimensional, and of minimal regular Hodge--Tate type, we construct an Euler system for $\rho$ in the sense of forthcoming work of Jetchev--\nek--Skinner. The construction is based on Theta cycles as introduced in a previous paper, following works of Kudla and Liu on arithmetic theta series on unitary Shimura varieties; it relies on a certain modularity hypothesis for those theta series. 

Under some ordinariness assumptions, one can attach to $\rho$ a $p$-adic $L$-function.  By recent results of Liu and the author, and the theory of Jetchev--\nek--Skinner,  we deduce the following (unconditional) result under mild assumptions:  if the $p$-adic $L$-function of $\rho$ vanishes to order $1$ at the centre, then the Selmer group of $\rho$  has rank~$1$, generated by the class of an algebraic cycle.  This confirms a case  of the $p$-adic Beilinson--Bloch--Kato conjecture.

\end{abstract}

\tableofcontents

 \section{Introduction}

A remarkable construction of Kolyvagin shows that if a Heegner point is non-torsion, then the Mordell--Weil and Selmer groups of a (modular) elliptic curve both have rank one \cite{Koly}. Combined with the formulas of Gross--Zagier and Perrin-Riou \cite{GZ, PR}, which relate heights of Heegner points and derivatives of $L$-functions in complex or $p$-adic coefficients, Kolyvagin's work provides important evidence for the Birch and Swinnerton-Dyer conjecture and its $p$-adic analogue.

We are interested in analogous pictures for higher-rank motives, or more simply  geometric\footnote{This and other unexplained notions will be defined in the main body of the paper.}  Galois representations, of weight~$-1$ . The most accessible ones  are arguably those over a CM field that are  \emph{conjugate-symplectic}. For those Galois representations, Jetchev--\nek--Skinner have recently theorised a  variant of Kolyvagin's method based on the notion of (what we propose to call) a \emph{JNS} Euler system; this is still a system of Selmer classes satisfying certain compatibility conditions. 
The purpose of this work is to construct such an Euler system, for those representations as above that are automorphic and even-dimensional of minimal regular Hodge--Tate type.

\medskip

 The companion  formulas of Gross--Zagier/Perrin-Riou type were recently proved in \cite{LL, LL2} and \cite{DL},  which  allows  to obtain various applications to the analogues (by Beilinson, Bloch, Kato, and Perrin-Riou) of the Birch and Swinnerton-Dyer conjecture.\footnote{In a future version, we plan to include applications to anticyclotomic 
Iwasawa theory.}

 In the rest of this introduction, we briefly state our main result and consequence, and the idea of its proof. For an overview  on the context and history of the constructions, and statements of other arithmetic applications, 
 we refer to  \cite{cet}.

\subsection{Main result}\lb{sec:ass}

Let $E$ be a CM field with absolute Galois group $G_{E}$ and maximal totally real subfield $F$, and let $\rc\in\Gal(E/F)$ be the complex conjugation. 

Let $n=2r$ be an even positive integer and let   
$$\rho\colon G_{E}\to \GL_{n}(\Qpb)$$ be  an  irreducible continuous representation,  that is geometric in the sense of \cite[I, \S1]{FM95}.
 We  denote by $\rho^{\rc}\colon G_{E}\to \GL_{n}(\Qpb)$ the representation defined by $\rho^{\rc}(g)=\rho(cgc^{-1})$, where $c\in G_{E}$ is any  fixed lift of~$\rc$. (A different choice of lift would yield an isomorphic representation.)

Suppose that the following conditions are satisfied:
\begin{enumerate}
\item \lb{csy}
$\rho $ is \emph{conjugate-symplectic} in the sense that there exists a perfect pairing
$$\rho\ot_{\Qpb}\rho^{\rc}\to \Qpb(1)$$
such that for the induced map $u\colon \rho^{c}\to \rho^{*}(1)$ (where ${}^{*}$ denotes the linear dual) and its conjugate-dual $u^{*}(1)^{\rc}\colon \rho^{\rc}\to \rho^{\rc, *}(1)^{\rc}=\rho^{*}(1)$, we have $u=-u^{*}(1)^{\rc}$;
\item\lb{HT}  for every place $w\vert p$ of $E$ and every embedding  $\jmath  \colon E_{w}  \into \mathbf{C}_{p}$, the $\jmath$-Hodge--Tate weights\footnote{Our convention is that  the cyclotomic character has weight~$-1$.}
of $\rho$ are the $n$ integers $\{-r, -r+1, \ldots,r-1\}$;
\item \lb{auto}
$\rho$ is \emph{automorphic} in the sense that for each $\iota\colon \Qpb \into \bC$, there is a cuspidal automorphic representation $\Pi^{\iota}$ of $\GL_{n}(\A_{E})$ such that $L_{\iota}(\rho, s)=L(\Pi^{\iota}, s+1/2)$;
\end{enumerate}
For a place $v$ of $F$, denote by $\rho_{v}$ the restriction of $\rho$ to $G_{E_{v}}:=\prod_{w\vert v}G_{E_{w}}$ (where the product ranges over the one or two places of $E$ above $v$).
For each ideal $m\subset \sO_{E}$, we have a ring class field $E[m]\supset E$; we also put $E[0]:=E$. We denote by $H^{1}_{f}(E[m], \rho)$ the Bloch--Kato Selmer groups \cite{BK}.

 Let $\sM_{1}$ be a  set consisting of all but finitely many of the places $v$ of $F$ that are split in $E$ and at which $\rho$ is unramified, and let $\sM$ be the set of finite subsets of $\sM_{1}$,  which we identify with a set of squarefree ideals in $\sO_{F}$.  Fix a set  $\wp$ of $p$-adic places of $F$ that are split in $E$, such that for each $v\in \wp$, the representation $\rho_{v}$ is Panchishkin-ordinary (Definition \ref{ord P gal}) and crystalline,\footnote{In fact the crystalline condition can be replaced by the considerably weaker condition of \eqref{cond: marcil} and even entirely removed, see Remark \ref{rema: marcil}.}   and let $\sM[\wp]$ be the set of ideals of the form $m\prod_{v\in \wp}v^{s_{v}}$ with $s=(s_{v})\in \Z_{\geq 0}^{\wp}$.

\begin{theoA}\lb{A} Let $\rho$ be a representation satisfying conditions 1., 2., 3. above, and let $\sM_{1}, \sM, \wp, \sM[\wp]$ be as above.  Assume that the root number $\eps(\rho)=-1$, that $F\neq\Q $  or $n=2$, and  that the Modularity Hypothesis \ref{mod} holds. 
 
The system of classes 
$$\Theta_{m}\in H^{1}_{f}(E[m], \rho), \qquad m\in \sM[\wp]\cup \{0\},$$
 of Definition \ref{def th c} forms a  JNS Euler system.
\end{theoA}
For the definition of  JNS Euler systems\footnote{Some of the (young) literature on the subject calls this notion  ``split anticylcotomic Euler system''.}
and the precise statement of the theorem, see Theorem \ref{norm rel}. 
\begin{rema}
Hypothesis \ref{mod} concerns the modularity of a certain generating series of Selmer classes coming from cycles in unitary Shimura varieties, for which the evidence is discussed in \cite[Remark 4.4]{cet} and references therein. We also  rely on a description of part of the cohomology of those varieties, Hypothesis \ref{hyp coh}, expected to be confirmed in a sequel to \cite{KSZ}.

The assumption on the root number is natural in the sense that, by \eqref{jns impl} below and the Beilinson--Bloch--Kato conjecture (e.g. \cite[Conjecture 2.2]{cet}), in the complementary case $\eps(\rho)=+1$ every JNS Euler system is expected to be zero.
\end{rema}

The main result of the work of Jetchev--\nek--Skinner (see \cite{Skinner} or \cite[\S~8]{ACR}) implies that, under mild conditions on the image of $\rho$, we have
\beq\lb{jns impl}
\Theta_{0} \neq 0 \qquad \Longrightarrow \qquad H^{1}_{f}(E, \rho)=\Qpb\Theta_{0}.\eeq
  Thus Theorem~\ref{A} demands a nonvanishing criterion for $\Theta_{0}$. Under some ramification restrictions:
\begin{itemize}
\item Li and Liu have proved a nonvanishing criterion in terms of derivatives of $L$-functions  \cite{LL, LL2},  conditionally on some standard conjectures on Abel--Jacobi map;
\item Liu and the author, under the  further assumption that one can take $\wp=\{\text{all places of $F$ above $p$}\}$,
have proved an unconditional nonvanisihing criterion in terms of $p$-adic $L$-functions, and confirmed the Modularity Hypothesis in that context \cite{DL}.
\end{itemize}
For the precise statements cast into the setup of the present paper, and their consequences towards the complex and $p$-adic Beilinson--Bloch--Kato conjectures,\footnote{See \cite[Conjecture 2.2]{cet} for a formulation.} see \cite[Theorem A]{cet}.  

We highlight  the following corollary also stated as part of \emph{loc. cit.}, which appears to be the first complete result towards the  Beilinson--Bloch--Kato conjectures in analytic rank~1 for high dimensional representations --  \emph{ex aequo} with a  similar result on Rankin--Selberg motives obtained  \cite[Theorem C]{DZ} (using \cite{LTXZZ} or \cite{Lai-Ski}, cf. Remark \ref{other wk}). It follows from combining \cite[Theorems 1.7, 1.8]{DL},  Theorem~\ref{A}  above, and the forthcoming theory of Jetchev--\nek--Skinner as in  \cite[Theorem 8.3]{ACR}, \cite{Skinner}. 

\begin{coro*}
 Suppose further that $E/F$ is totally split above $2$ and  $p$, that $p>n$, that places of $F$ ramified in $E$ are unramified over $\Q$,  and that the representation $\rho$ is:
 \begin{itemize}
 \item Panchishkin-ordinary and crystalline at all $p$-adic places of $E$;
 \item of `large image' in the sense that it is absolutely irreducible and there exists a $\gamma\in\Gal(\overline{E}/E(1)(\mu_{p^{\infty}}))$ such that $\dim_{\Qpb} \rho /(\gamma-1) \rho=1$;
 \item
 `mildly ramified'   in the sense that the associated  automorphic representation  $\pi$  (\S~3.2) satisfies \cite[Assumption 1.6 (1)-(2)-(3)]{DL}.
 \end{itemize}
 Denote by $\sX_{F}$ the $\Qpb$-scheme of continuous $p$-adic characters of $G_{F}$ that are unramified outside $p$, by $\frak m\subset\sO(\sX_{F})$ the ideal of functions vanishing at $\one$, and by $L_{p}(\rho)\in \sO(\sX_{F})$  the $p$-adic $L$-function of $\rho$ from \cite{DL} (see \cite[Proposition 5.2]{cet}). 
 
 Then 
$$\textup{ord}_{\frakm} L_{p}(\rho)=1 \quad \Longrightarrow\quad  \dim_{\Qpb}H^{1}_{f}(E, \rho)=1.$$
\end{coro*}

\begin{enonce}[remark]{Example} Suppose that $A$ is a modular elliptic curve over $F$ without complex multiplication, and let $n=2r$ be an even positive integer.
Suppose that $E/F$, $p$, and $n$ satisfy the conditions of the corollary, and that:
\begin{itemize}
\item $G_{F}$ surjects onto ${\r{Aut}}_{\mathbf{F}_{p}}A[p]$ (which happens for every sufficiently large prime $p$ by Serre's open image theorem);
\item $A$ has good ordinary reduction at all places above $p$ (which happens for a density-one set of primes $p$);
\item all places of bad reduction of $A$ split in $E$.
\end{itemize}
Then 
$$\textup{ord}_{\frakm} L_{p}({\r{Sym}}^{n-1}V_{p}A_{E}(1-r))=1 \quad \Longrightarrow\quad  \dim_{\Q_{p}}H^{1}_{f}(E, {\r{ Sym}}^{n-1}V_{p}A_{E}(1-r))=1.$$
Indeed, $\rho:={\r{Sym}}^{n-1}V_{p}A_{E}(1-r))$ verifies the conditions of Panchishkin-ordinary and crystalline and `mildly ramified' by \cite[\S~1.4]{DL}. For the condition of `large image', the absolute irreducibility condition follows from the irreducibility of symmetric power representations of $\GL_{2}$. It remains to find an element $\gamma$ as required in the corollary. Let $\gamma\in G_{F}$ be an element acting on $T_{p} A $ by a nontrivial unipotent matrix $\rho_{A}(\gamma)$. Then  the action of $\gamma$ on $\mu_{p^{\infty}}$ is given by $\det\rho_{A}(\gamma)=1$. Up to replacing $\gamma$ by a power, we may assume it fixes $E[1]$ as well. It is clear that $\gamma-1$ has $1$-dimensional cokernel on  $\rho$.
\end{enonce}

\begin{rema} \lb{other wk} To the author's knowledge, the vast literature  on Euler systems contains only four other constructions  for high rank motives:  one by Liu--Tian--Xiao--Zhang--Zhu \cite{LTXZZ} for conjugate-symplectic  Rankin--Selberg motives, which is of a type introduced by Bertolini--Darmon in \cite{BD}; one by Cornut \cite{cornut}, for base-changes of some symplectic motives, of a type  similar to the one of \cite{Koly}; and two of JNS type, valid for an infinite range of Hodge--Tate weights: one  by Graham--Shah \cite{g-s}, for conjugate-symplectic motives that are also symplectic, and one (contemporaneous to the present paper) by Lai--Skinner \cite{Lai-Ski}, again for conjugate-symplectic  Rankin--Selberg motives.
\end{rema}

 \subsection{Idea of the proof}
 The construction of the  \emph{Euler system of Theta cycles} starts from the arithmetic theta lifts on unitary Shimura varieties introduced by   Liu in \cite{Liu11} (partly based on a construction of  Kudla \cite{Kud97}); as in \cite{cet}, we recast them as trilinear forms valued in $H^{1}_{f}(E, \rho)$. The higher layers of the system are given by taking connected components of the special cycles arising in the constructions, and varying the input data in a well-chosen way. 
 
 To prove that the constructed classes indeed form an Euler system we need to establish that they are integral and that they are bound up by  certain norm  relations (`horizontal' and `vertical', i.e. at non-$p$-adic and $p$-adic places). Following an idea pioneered in \cite{YZZ} and developed in the context of Euler systems in  \cite{LSZ}, we prove the horizontal norm relations based on the fact that the space of (scalar-valued) trilinear forms appearing in the constructions decomposes into a product of local spaces, each of dimension~1 by the theory of the local theta correspondence. Then some equivalent relations may be established in any models of these local spaces: in our context, an explicit one is given by the zeta integrals used by Godement--Jacquet to construct the standard $L$-functions for $\GL_{n}$, where  the desired identity  is easy to prove. This model in fact guides the choice of input data away from $p$; the integrality relations are then established by explicit computation.
      At $p$-adic places, we use a variant of choices of data from \cite{DL}, and prove its local nontriviality again by a computation in the Godement--Jacquet model.

 \medskip
 
In \S~\ref{sec 2}, we construct the system and reduce its fu-ndamental properties to local statements. In \S~\ref{sec 3}, we prove those statements.

\subsection*{Acknowledgements}
I would like to thank  Francesc Castella, Henri Darmon, Nadya Gurevich,  Dimitar Jetchev,  David Marcil, Christopher Skinner, Ariel Weiss, Wei Zhang, and especially Yifeng Liu and   Waqar Ali Shah  for useful conversations or correspondence. 

A substantial part of this paper was written during the KUMA International Summer School in Sarajevo in August 2021, and I am grateful to  its director Claudia Zini and all the staff for their hospitality.

\section{The Euler system of Theta cycles}
\lb{sec 2}
\subsection{Setup}
We briefly review the setup for the construction of Theta cycles, referring to \cite{cet} and references therein for the details. 

\subsubsection{Notation}
Suppose for the rest of this paper that $E$ is a CM field with totally real subfield $F$.   We denote by ${\r c}\in \Gal(E/F)$  the complex conjugation, and by $\eta\colon F^{\ts}\bks\A^{\ts}\to\{\pm 1\}$   the quadratic character attached to $E/F$. 
 We denote by $\A$ the ad\`eles of $F$; if $S$ is a finite set of places of $F$, we denote by $\A^{S}$  the ad\`eles of $F$ away from $S$. If $\G$ is a group over $F$ and $v$ is a place of $F$, we write $G_{v}\coloneqq \G(F_{v})$; if $S$ a finite set of places of $F$, we write $G_{S}\coloneqq \prod_{v\in S}G(F_{S})$. (For notational purposes, we will identify a place of $\Q$ with the set of places of $F$ above it.) We denote by  $\psi\colon F\bks \A\to \bC^{\ts}$ the standard additive character with $\psi_{\infty}(x)=e^{2\pi i \Tr_{F_{\infty}/\R}x}$,  and we set $\psi_{E}\coloneqq \psi\circ \Tr_{E/F}$. 
 
We fix a rational prime $p$ and  denote by $\Q^{\cyc}\subset \Qpb$ the extension of $\Q$ generated by all roots of unity. 

  We fix an embedding $\iota^{\cyc}\colon \Q^{\cyc}\into \bC$,  by which we view $\psi_{|\A^{\infty}}$ as valued in $\Q^{\cyc}$.  We denote by $\sO$ the integral closure of $\Z_{p}$ in $\Qpb$.

\subsubsection{Quasisplit unitary group}
Let $W=E^{n}=W^{+}\oplus W^{-}$  where $W^{+}=\Span(e_{1}, \ldots, e_{r})$, $W^{-}=\Span(e_{r+1}, \ldots, e_{2r})$,  equipped with  the skew-hermitian form $\lan\ , \ \ran_{W}$ with matrix  $\smalltwomat{}{ 1_{r}}{-1_{r}}{}$ (here $1_{r}$ is the identity matrix of size $r$). We denote by $\G=\rU(W)$ its unitary group, which we may view as a subgroup of $\Res_{E/F}\GL_{n}$. Denote by $\rP\subset \G$ the parabolic subgroup stabilizing $W^{-}$, and by $\Herm_{r}$ the space of hermitian $r\ts r$ matrices.

We have:
     \begin{itemize}
     \item a Weyl element 
     $$w=\twomat{}{1_{r}}{-1_{r}}{}\in \G;$$
       \item a homomorphism 
\beqq
m\colon\Res_{E/F}\GL_r&\to \rP\subset \G\\
a&\mapsto          
          m(a)\coloneqq
          \begin{pmatrix}
              a &  \\
               & {}^{\rt}{a}^{\rc,-1}. \\
          \end{pmatrix}
\eeqq
          whose image is a Levi factor of $\rP$;
       \item a homomorphism 
       \beqq
       n\colon\Herm &\to \rP\subset\G\\
b&\mapsto          
          n(b)\coloneqq
          \begin{pmatrix}
              1_r & b \\
               & 1_r \\
          \end{pmatrix}
          ,
\eeqq
whose image is the unipotent radical of $\rP$. Here, we denote by $\Herm$ the space of hermitian matrices; we will also denote by $\Herm(F)^{+}\subset \Herm(F)$ the subspace consisting of totally positive semidefinite matrices. 
     \end{itemize}

Attached to $\G$, we have:
\begin{itemize}
\item a $\Qpb$-vector space $\sH_{\Qpb}$ of  modular forms  (see \cite[\S~2.2]{DL}, \cite[\S~4.3]{cet}); 
\item for any ring $R$, the space $\ul{\r{SF}}_{\, R}$ of those formal (Siegel--Fourier) expansions
$$\sum_{T\in \Herm_{r}(F)^{+}} c_{T}(a) \, q^{T}, \qquad c_{T}\in C^{\infty}( \GL_{r}(\A^{\infty}_{E}) , R)$$
satisfying $c_{{}^{\rt}a^{\rc}Ta}(y)=
c_{T}(ay)$ for all $a\in \GL_{r}(E)$;
\item  an injective $p$-adic $q$-expansion map
\beq\lb{ul q}
\ul{\qqq}\colon \sH_{\Qpb} \longrightarrow \ul{\r{SF}}_{\, \Qpb},
\eeq
denoted by $\underline{\bf q}_{p}$ in \cite[\S~4.2]{cet}.
\end{itemize}

\subsubsection{Incoherent unitary groups}
Let $V$ be an incoherent, totally positive definite $E/F$-hermitian space;
  this is simply a collection of $E_{v}/F_{v}$-hermitian spaces $V_{v}$ of the same dimension, indexed by the places of $F_{v}$, which is not isomorphic to one of the form  $(V_{0}\ot_{F}F_{v})_{v}$ for some hermitian space $V_{0}$ over $E$, and such that $V_{v}$ is positive definite for all $v\vert \infty$.  
If $S$ is a finite set of places of $F$, we  put $V_{\A^{S}}=\ot_{v\notin S}V_{v}$. For $x_{1}, \ldots , x_{r}\in V_{v}$, we  have the moment matrix 
  $$T(x):=((x_{i}, x_{j})_{V_{v}})_{ij}\in \Herm_{r}(F_{v}).$$ 
  We denote  by $\H_{V}$ the  incoherent unitary group associated with $V$  in the sense of \cite{cet}. 

Attached to $\H_{V}$ we have a tower of  Shimura varieties 
$$(X_{\H_{V},K})_{K\subset \H_{V}(\A^{\infty})}$$
 of dimension $\dim V-1 $ over~$E$ as in \cite[\S~4.2]{cet}.

\subsubsection{Weil representation}\lb{weil} Let $v$ be a finite place of $F$, and let $V_{v}$ be an $E_{v}/F_{v}$-hermitian space of dimension~$n$. 
The basis $\{e_{1}, \ldots , e_{r}\}$ of $W^{+}$ identifies $V_{v}\ot_{E_{v}}W_{v}^{+}=V_{v}^{r}$. We have a representation $\omega_{v}:=\omega_{V_{v}}$ of $G_{v}\ts H_{V_{v}}$ on  the Schwartz space $\sS(V_v\ot_{E_{v}} W_{v}^{+},\Qpb)$,  characterized by the property that for  $\phi\in \sS(V_v\ot_{E_{v}}W_{v}^{+},\Qpb)$:
   \begin{itemize}
           \item for $h\in H_{V_{v}}$, we have
           \[
           \omega_{v}(h)\phi(x)=\phi(h^{-1}x);
           \]
         \item for $a\in \GL_{r}(E_{v})$ and $b\in\Herm_r(F_v)$, we have
\beqq
           \omega_{v}(m(a))\phi(x)&=|\det a|_E^r\cdot\phi(x a), \\
           \omega_{v}(n(b))\phi(x)&=\psi_{v}(\r{Tr} b T(x))\phi(x),\\
           \omega_{v}\(w_{r}\)\phi(x)&=\gamma_{V_v,\psi_{v}}^r\cdot\widehat\phi(x),\eeqq
           where            $\gamma_{V_v,\psi_{v}}\in\{\pm 1\}$ is the Weil constant of $V_v$ with respect to $\psi_{v}$,
           and $\widehat\phi$ denotes the Fourier transform
$$            \widehat\phi(x)\coloneqq\int_{V_v^r}\phi(y)\psi_{E,v}\(\sum_{i=1}^r(x_i,y_i)_V\) \, \rd y
$$
for the $\psi_{E,v}$-self-dual Haar measure $\rd y $ on $V_v^r$.
      \end{itemize}

For an incoherent $E/F$-hermitian space $V$, we  put $\omega=\omega_{V}=\ot_{v}\omega_{V_{v}}$, the product running over all finite places of $F$; it is a representation of $\G(\A^{\infty})\ts \H_{V}(\A^{\infty})$ on the space $\sS(V_{\A^{\infty}}\ot_{E}W^{+}, \Qpb)$.

 \subsubsection{$p$-adic automorphic representations} \lb{p aut rep}
Given our Galois representation $\rho$, we choose a relevant $p$-adic automoprhic representation of $\G(\A)$  over $\Qpb$ (in the sense of \cite[Definition 3.2]{cet})
$$\pi\subset \sH_{\Qpb}$$
whose base-change to $\GL_{n}(\A_{E})$ is $\Pi_{\rho}$, which exists by \cite[Proposition 3.4]{cet}. (The representation $\pi$ is not uniquely determined by this condition, although as noted in \cite[Remark 3.5]{cet}, there is a `standard' choice.)

We enforce from now on the assumption that $\eps(\rho)=-1$. Then by \cite[Proposition 3.8]{cet}, attached to $\rho$ and $\pi$ we have  a pair consisting of 
 \begin{itemize}
 \item  an incoherent   totally definite $E/F$-hermitian space $V$ of dimension $n$, and 
\item a relevant $p$-adic automorphic representation $\sg$ of $\H(\A)$ over $\Qpb$  (in the sense of \cite[Definition 3.2]{cet})
\end{itemize}
uniquely characterised (up to isomorphism) by the condition that  the space of coinvariants
$$\Lm_{\rho}\coloneqq  (\pi^{\vee}\ot\omega\ot\sg)_{\G(\A^{\infty})\ts \H_{V}(\A^{\infty})} $$
is nonzero. The space $\Lm_{\rho}$  is then in fact $1$-dimensional over $\Qpb$. 

Henceforth, we write $\H=\H_{V}$, $\omega=\omega_{V}$,  $X=X_{\H_{V}}$.

  \subsubsection{Realisation of $\sg$ in cohomology} 
 From now on we assume that $F\neq \Q$ or $n=2$, which implies that the varieties $X_{\H_{V}, K}$ are projective -- except in a case related to modular curves where $X_{\H_{V}, K}$ can be canonically compactified  by adding finitely many cusps; in  that case, we replace $X_{\H_{V}, K}$ by its compactification.

For each open compact $K\subset \H(\A^{\infty})$, let 
$$\sg_{\rho}^{K}:= \Hom_{\Qpb[G_{E}]} (H^{2r-1}_{\textup{\'et}}(X_{\H_{V}, K,\ol E}, \Qpb(r)), \rho).$$
 We will assume the following hypothesis (a variant of \cite[Hypothesis 4.1]{cet}); it is known for $n=2$, and  it is expected to be confirmed in general in a sequel to \cite{KSZ}.
\begin{enonce}{Hypothesis}   \lb{hyp coh}
For each open compact subgroup $K\subset \H(\A^{\infty})$,  we have an isomorphism of $\Qpb[K\bks \H_{V}(\A^{\infty})/K]$-modules
\beq \lb{eq coh}
\sg_{\rho}^{K}\cong \bigoplus_{\sg'} \sg'^{K},
\eeq
where the direct sum runs over the isomorphism classes of  relevant $p$-adic automorphic representation (in the sense of \cite[Definition 3.2]{cet})
 $\sg'$ of $\H_{V}(\A)$ with ${\r{BC}}(\sg')=\Pi$.
\end{enonce}
We put $M_{\sg, K}:= \sg^{\vee,K}\ot \rho$, which we identify with a subspace of  $M_{\rho, K}:= \sg_{\rho}^{\vee, K}\ot\rho\subset H^{2r-1}_{\textup{\'et}}(X_{K,\ol E}, \Qpb(r)).$ 
We then identify  
$$
\sg= \varinjlim_{K} \Hom_{\Qpb[K\bks \H_{V}(\A^{\infty})/K]}(M_{\sg, K}, \rho)\subset  \varinjlim_{K} \Hom_{\Qpb[K\bks \H_{V}(\A^{\infty})/K]}(M_{\rho, K}, \rho) = \varinjlim_{K} \sg_{\rho, K}.$$

\Subsection{Special cycles and generating series over ring class fields}
\subsubsection{Connected components of unitary Shimura varieties}
Let $\rT$ be the unitary group of $E$ with the form induced by the norm $N_{E/F}$ (as an algebraic group over $F$), and denote  $$\sC:=\{ \text{open compact subgroups } C\subset \rT(\A^{\infty})\}\cup \{\rT(\A^{\infty})\}.$$ (Elements $C\in \sC$ will often appear as sub/superscript, omitted when $C=\rT(\A^{\infty})$.) For each  $C\in \sC$,  let $E_{C}$ be the abelian extension of $E$ with 
$$\Gal(E_{C}/E) =: \Gamma_{/C}:= \rT(F)\bks \rT(\A^{\infty})/C$$
 under the class field theory isomorphism (which we will view as an identification). Let
  $\Gamma:= \varprojlim_{C}\Gamma_{/C}$.
  For any profinite group $\Gamma'$, we will denote by $\widehat{\Gamma}:=\varinjlim_{C'}\Spec\Qpb[\Gamma'/C']$ the space of locally constant $\Qpb$-valued characters of $\Gamma'$ (where the limit ranges over finite-index subgroups).

\medskip

Let $X_{\rT}$ be the  tower of $0$-dimensional Shimura varieties over $E$ associated with .
 For every coherent or incoherent unitary group $\H'$, 
we denote by $\nu_{\rH'}\colon \rH'\to \rT$ the determinant character (the subscript will be omitted when understood from the context).  
The tower $X_{\rH'}=(X_{\rH', K'})_{K'}$ of Shimura varieties maps to the tower $X_{\rT}$ via surjective morphisms still denoted
 $$\nu\colon X_{\rH, K'}\to X_{\rT, \nu(K')}.$$
These induce  bijections on the set of geometrically connected components.

Fix an identification of  $\Gamma$-sets  $X_{\rT}(E^{\mathrm{ab}})\cong \Gamma$. A subset $S\subset \rT(\A^{\infty})$ is said to be of \emph{level} $C\in\sC$ if $C$ is minimal for the property that $S$ is a union of $C$-cosets. 
If  $S\subset \rT(\A^{\infty})$ is of level  $C\supset C'$, let
$$X_{\rT,C'}^{S}\subset X_{\rT, C', E_{C}}$$
be the  $E_{C}$-subscheme whose set of $E^{\mathrm{ab}}$-points is identified with the image of  $S$ in $\Gamma_{/C'}$. If  $S\subset \rT(\A^{\infty})$ is of level   $C\supset \nu(K)$, let 
$$X_{\H, K}^{S}:=\nu^{-1}(X_{\rT, \nu(K)}^{(S)}) \subset X_{\H, K, E_{C}}.$$

Then for each $i,j$,  each $C\in \sC$, and each $t\in \Gamma_{/C}$ we have a  direct $\Q_{p}[G_{E_{C}}]$-module summand
$H^{i}_{\text{\'et}}(X^{tC}_{K,\ol{E}}, \Q_{p}(j))\subset 
   H^{i}_{\text{\'et}}(X_{K,\ol{E}}, \Q_{p}(j)), $
   and we denote by
\beq \lb{dir sum}
\rr^{tC}\colon   H^{i}_{\text{\'et}}(X_{K,\ol{E}}, \Q_{p}(j)) \to    H^{i}_{\text{\'et}}(X^{tC}_{K, \ol{E}}, \Q_{p}(j))\eeq
the projection induced by the inclusion $X_{K}^{tC}\into X_{ K, E_{C}}$.

 \subsubsection{Special cycles} \lb{sec: sp cy}
From now on, we abbreviate $\rH=\rH_{V}$ and $X_{?}:=X_{\H, ?}=X_{\H_{V}, ?}$ (for any decoration `?').
For a $C\in \sC$ and a compact open subgroup $K'\subset \H'(\A^{\infty})$ (for some unitary group $\H'\stackrel{\nu}{\to}\rT$), we write $K'^{C}:= K'\cap \nu^{-1}(C)$.

Let $x\in V_{\A^{\infty}}\ot_{E} W^{+} =V_{\A^{\infty}}^{r}$.
\begin{itemize}
\item
Suppose that
\beq\lb{x adm}
T(x):=((x_{i}, x_{j})_{V})_{ij}\in \Herm_{r}(F)^{+} \text{\quad and \quad} V(x):=\Span_{E}(x_{1}, \ldots, x_{r}) \text{ is positive-definite}.
\eeq
For any compact open subgroup $K\subset \G(\A^{\infty})$, any $C\in \sC$, and any $t\in \rT(\A^{\infty})$, 
we have a cycle $Z^{tC}(x)_{K}$ defined as follows (see  \cite[\S~3A]{Liu11} or \cite[\S~4]{LL} for more details when $tC=\rT(\A^{\infty})$). 

Pick an embedding $\iota\colon E\into \bC$, and let $V^{\iota}$ be the (unique up to isomorphism) totally definite hermitian space over $E$ with $V^{\iota}_{\A^{\infty}}\cong  V_{\A^{\infty}}$; we fix such an identification. Then we may write $x=h^{-1}x'$ for $x'\in (V^{\iota})^{r}$ and $h\in \H(\A^{\infty})$. 
 Let $\H(x')$ be   the  unitary group of  the  subspace $V(x')^{\perp}\subset V^{\iota}$, where $V(x')
 :=\Span_{E}(x'_{1}, \ldots, x_{r}') $;
let $K_{x'} := hKh^{-1}\cap \H(x')(\A^{\infty})$.  We also denote by $\H(x)(\A)$  the unitary group of $V(x)^{\perp}$, and set $K_{x}:=K\cap \H(x)(\A)$.

The natural inclusion 
$j_{x'}\colon \H(x'){\into} \H^{\iota}$ of unitary groups induces a morphism
  of Shimura varieties 
 $$ X_{\H(x'), K_{x'}}\stackrel{j_{x', K}}{\longrightarrow} X_{hKh^{-1}}\stackrel{\cdot h}{\longrightarrow}X_{K}.$$
Let $\sL_{K}$ be the Hodge bundle on 
$X_{K}$ (or its base-change to $E_{C}$).
We  then define a cycle (see \cite[\S~2.5]{Ful} for general background)
 $$Z^{tC}(x)_{K} :=
   c_{1}(\sL^{\vee}_{K})^{\dim V(x)-r}_{|X^{tC}_{K}}
  \frown  [C\nu(K_{x}) : C]^{-1}\cd 
[j_{x, K}(X_{\rH(x), K_{x}}^{tC\nu(K_{x})}) ]
 \in \Ch_{r-1}(X_{K, E_{C}})_{\Q}, $$ 
 where we have used the suggestive notation\footnote{This notation can likely be given a substantive meaning in the framework of \cite{ST}.}
\beq\lb{sugg}
j_{x, K}(X_{\rH(x), K_{x}}^{S})
:=
 j_{x', K}(X_{\rH(x'), K_{x'}}^{\nu^{-1}(h) S}) h, \eeq
in which the right multiplication denotes the action of $\H(\A^{\infty})$ on the tower $(X_{\H, K})_{K}$. 

The definition is independent of the auxiliary choices made.

\item   If  $x$ does not satisfy \eqref{x adm},
we put $  Z^{C}(x)_{K}:=0$. 
\end{itemize}

 It is clear that $\Tr_{E_{C'}/E_{C}}   Z^{C'}(x)_{K} = Z^{C}(x)_{K}$ whenever $C'\subset C\in \sC$.

For a locally constant function $\chi\colon  { \Gamma}\to \Qpb$, we also define 
$$Z(x, \chi)_{K}:= \sum_{t\in \Gamma_{/C}}\chi(t) Z^{tC}(x)_{K}
\qquad
 \in \Ch_{r-1}(X_{K, E_{C}})_{\Qpb},$$
for any $C\in \sC$ such that $\chi$ factors through $\Gamma_{/C}$. (Thus $Z^{C}(x)_{K}=Z(x, \one_{C})_{K}$ for any $C\in \sC$.)

\begin{rema}
Suppose that $\chi\colon  { \Gamma}\to \Qpb^{\ts}$ is a locally constant character. For every  $x\in V_{\A^{\infty}}^{r}$, $\gamma\in \Gamma_{}$, and $h\in \H(x)(\A^{\infty})$, we have
\beq\lb{chi act tht}              
Z(x, \chi)_{K}^{\gamma}=\chi^{-1}(\gamma)Z(x, \chi)_{K}.\eeq
and 
\beq\lb{Zchi eq}
Z^{C}(x)_{K}h = Z^{\nu(h)C}(h^{-1}x)_{h^{-1}Kh}, \qquad Z(x,\chi)_{K}h = \chi^{-1}_{\H}(h) Z(h^{-1}x,\chi)_{h^{-1}Kh}.
\eeq
where we still denote simply as a right multiplication the pushforward action on cycles induced by the right action of $\H(\A^{\infty})$ on  $(X_{\H, K})_{K}$. 
\end{rema}

\subsubsection{Projection to the $\rho$-component}
Let $C\in \sC$, and let $K\subset \H(\A^{\infty})\cap \nu^{-1}(C)$ be an open compact subgroup. 
Denote by ${\r{Fil}}^{\bullet, C}\subset H^{2r}_{\textup{\'et}}(X_{ K, E_{C}}, \Q_{p}(r)) $ the filtration induced by the Hochschild--Serre spectral sequence $H^{i}(E_{C}, H^{2r-i}_{\text{\'et}}(X_{ K, E_{C}}, \Q_{p}(r))) \Rightarrow  H^{2r}_{\text{\'et}}(X_{ K, E_{C}}, \Q_{p}(r)) $. 
We have an absolute cycle class map
$$\mathrm{AJ} \colon \Ch_{r-1}(X_{K, E_{C}})_{\Qpb} \to  H^{2r}_{\textup{\'et}}(X_{K, E_{C}}, \Q_{p}(r)) / {\r{Fil}}^{2, C}.$$

\begin{lemm} \lb{p rho}
The Hecke-eigenprojection 
$$e_{\rho}\colon \bigoplus_{i\in \Z} H^{i}_{\textup{\'et}}(X_{ K, \ol E}, \Q_{p}(r)) \to M_{\rho, K}$$
 induces a Hecke-equivariant projection, still denoted 
$$ e_{\rho}\colon  H^{2r}_{\textup{\'et}}(X_{K, E_{C}}, \Q_{p}(r)) /\r{Fil}^{2, C} \to H^{1}(E_{C},  M_{\rho, K}),$$
such that the composition
\begin{gather*}
(-)_{\rho}\colon \Ch_{r-1}(X_{K, E_{C}})_{\Q} \stackrel{\mathrm{AJ}}{\longrightarrow} 
  H^{2r}_{\textup{\'et}}(X_{K , E_{C}}, \Q_{p}(r)) / {\r{Fil}}^{2, C}
\stackrel{e_{\rho}}{\longrightarrow} 
H^{1}(E_{C},  M_{\rho, K})\\
Z\mapsto Z_{\rho}:= e_{\rho}\mathrm{AJ}(Z)
\end{gather*}
takes values in $H^{1}_{f}(E_{C},  M_{\rho, K})$.
\end{lemm}

\begin{proof} 
As in \cite[Lemma 4.2]{cet}.
\end{proof}

  \subsubsection{Generating series} 
  Let $\chi\colon \Gamma \to \Qpb$ be a locally constant function, and let $C(\chi)\in\sC$ be maximal such that $\chi$ factors through $\Gamma_{/C(\chi)}$. 
For $\phi\in \sS(V_{\A^{\infty}}^{r})$ and any $K\subset \H(\A^{\infty})$ fixing $\phi$, let
$$Z_{T}(\phi, \chi)_{K}:= \sum_{x\in K\bks V_{\A^{\infty}}^{r} \ : \ T(x)=T} 
\phi(x) Z(x, \chi)_{K} \quad  \in \Ch_{r-1}(X^{}_{K})_{\Q};$$ 
 in the special case  $\chi=\one_{tC}$ for $C\in \sC$, we write 
$$Z^{tC}_{T}(\phi):=Z_{T}(\phi, \one_{tC}), \qquad {}^{\qqq} \Theta^{tC}(\phi)_{\rho, K}= {}^{\qqq} \Theta(\phi, \one_{tC})_{\rho, K}.$$

We define
\beqq
{}^{\qqq} \Theta(\phi, \chi)_{\rho, K}(a)
&:= {\vol(K)}\sum_{t\in \Gamma/{C}} \chi(t)
\sum_{x\in K\bks V_{\A^{\infty}}^{r}} \phi(xa) Z^{\nu(m(a)) t C}(x)_{K, \rho} \, q^{T(x)}  
&\in 
  H^{1}_{f}(E_{C(\chi)}, M_{\rho, K})
\ot_{\Q_{p}}{\underline{\mathrm{SF}}}_{\Qpb},
\eeqq 
where $\vol$ is as in \cite[Definition 3.8]{LL}.

\begin{lemm}\lb{pushfwd}
Let $K'\subset K\subset \H(\A^{\infty})$ be compact open subgroups,
 let
  $\phi\in \sS(V_{\A^{\infty}}^{r})^{K}$, and let  $\rp=\rp_{K'/K}\colon X_{K'}\to X_{K}$ denote the projection map. 
Then for every locally constant $\chi\colon \Gamma \to \Qpb$, we have
$$\vol(K')\,  \rp_{*} Z_{T}(\phi, \chi)_{K'} = \vol(K) \,   Z^{}_{T}(\phi, \chi)_{K}.$$
\end{lemm}
\begin{proof}
It suffices to show that   $\rp_{*} Z^{C}_{T}(\phi)_{K'}  =[K:K']\,  Z^{C}_{T}(\phi)_{K} $ for every $C\in \sC$.

Suppose first that, in the notation of  \S\ref{sec: sp cy}, we may identify $x$ with an element of $(V^{\iota})^{r}$. The finite surjective map $X_{\rH(x), K'_{x}}^{C\nu(K'_{x})} \to  X_{\rH(x), K_{x}}^{C\nu(K_{x})}$  has degree  
$
 [K_{x}:K_{x}'] \cdot [C\nu(K_{x}) : C\nu(K_{x}')]^{-1}. $
Since
$$
\rp^{*}\sL_{K}^{\vee}= \sL_{K'}^{\vee},$$
by the definitions and the projection formula,  we find  
$$ \rp_{*} Z^{C}(x)_{K'}
 =   [K_{x}:K_{x}'] \cdot[C\nu(K_{x}) : C\nu(K_{x}')]^{-1} {[C\nu(K_{x}'):C]^{-1}} \, [ j_{x, K}(X_{\rH(x), K_{x}}^{C\nu(K_{x})})] 
= [K_{x}:K_{x}'] \,   Z^{C}(x)_{K}.$$
It is easy to verify that this result remains valid without the assumption $x\in (V^{\iota})^{r}$. 

 Then, setting \beq\lb{phiT}\phi_{|T}(x) := \one_{[T(x)=T]} \phi(x),\eeq we find
\beqq
  \rp_{*} Z^{C}_{T}(\phi)_{K'}  = 
 \sum_{x\in K'\bks V_{\A^{\infty}}^{r} }   \phi_{|T}(kx)  \rp_{*} Z^{C}(kx)_{K'}
=
  \sum_{x\in K\bks V_{\A^{\infty}}^{r} }\sum_{k\in K'\bks K/ (K_{x}'\bks K_{x})}    \phi_{|T}(kx)    [K_{kx}:K_{kx}']    Z^{C}(kx)_{K} .
\eeqq
Now all the last three terms are independent of $k$, so that the inner sum equals
$$\sum_{x\in K\bks V_{\A^{\infty}}^{r} }     [K:K']   \phi_{|T}(x) 
   Z^{C}(x)_{K} = [K:K']  Z_{T}^{C}(\phi)_{K},
$$
as desired.
\end{proof}

\begin{coro}\lb{cor pushfwd}
The construction of ${}^{\qqq} \Theta(\phi, \chi)_{\rho, K}$ is compatible under pushforward in the tower $(X_{ K})_{K}$.
\end{coro}

\subsubsection{Modularity}
For the history and evidence in favour of the following  conjecture (which is \cite[Conjecture 4.17]{DL}),\footnote{The formulation in \emph{loc. cit.} is slightly different but easily seen to be equivalent.} see   \cite[Remark 4.4]{cet} and references therein.

\begin{enonce}{Hypothesis}[Modularity]  \lb{mod}
For every $\phi\in   \sS(V_{\A^{\infty}}^{r})$ and every $K\subset \H(\A^{\infty})$ fixing $\phi$, there exists a unique 
$$\Theta(\phi)_{\rho, K} \in 
H^{1}_{f}(E,  M_{\rho, K}) \ot_{\Qpb} \sH_{\Qpb^{}}$$
such that
$$\underline{\qqq} (\Theta(\phi)_{\rho, K}) =  {}^{\qqq}\Theta(\phi)_{\rho, K}.$$
\end{enonce}

We can amplify the modularity  to the other generating series.

\begin{prop} Assume Hypothesis \ref{mod}. Then for every locally constant function
$\chi\colon {\Gamma}\to \Qpb$,   every $\phi \in \sS(V_{\A^{\infty}}^{r})$, and every $K\subset \H(\A^{\infty})$ fixing $\phi$, there exists a unique 
$$
\Theta(\phi, \chi)_{\rho, K} \in 
H^{1}_{f}(E_{\chi},  M_{\rho, K})
\ot_{\Qpb} \sH_{\Qpb^{}}$$
such that 
$$
\underline{\qqq}(\Theta(\phi, \chi)_{\rho, K} )=  {}^{\qqq}\Theta(\phi, \chi)_{\rho, K}.
$$
\end{prop}
As usual, we will write $\Theta^{tC}(\phi)_{\rho, K} := \Theta(\phi, \one_{tC})_{\rho, K}$.

\begin{proof}
By Corollary \ref{cor pushfwd}, we may assume   that  $K$ satisfies $\nu(K)\subset C$.  
We define a slightly different  Siegel--Fourier expansion by 
$${}^{\qqq} \wtil{\Theta}(\phi, \chi)_{\rho, K}(a)
:= {\vol(K)}
\sum_{x\in K\bks V_{\A^{\infty}}^{r}} 
\omega(m(a))
\phi(xa) Z(x, \chi)_{K, \rho} \, q^{T(x)}  $$ 
and we put $ {}^{\qqq} \wtil{\Theta}^{tC}(\phi)_{\rho, K}={}^{\qqq} \wtil{\Theta}(\phi, \one_{tC})_{\rho, K}$.
It suffices to prove the proposition when $\chi$ is a character, in which case 
$${}^{\qqq} {\Theta}(\phi, \chi)_{\rho, K}(a)
=\chi(\nu_{\G}(m(a)))^{-1} \wtil{\Theta}(\phi, \chi)_{\rho, K}(a),$$
so that if $ \wtil{\Theta}(\phi, \chi)_{\rho, K}$ is a (Selmer-group-valued) Siegel modular form with $q$-expansion ${}^{\qqq} \wtil{\Theta}(\phi, \chi)_{\rho, K}$, then
 $$ {\Theta}(\phi, \chi)_{\rho, K}(g):= \chi^{-1}\circ\nu_{\G}\ot \wtil{\Theta}(\phi, \chi)_{\rho, K}$$ is a Siegel modular form with $q$-expansion ${}^{\qqq} {\Theta}(\phi, \chi)_{\rho, K}$.
Therefore it is equivalent to prove the modularity of the series ${}^{\qqq} \wtil{\Theta}(\phi, \chi)_{\rho, K}$ for all $\chi$, and  we may restrict to   $\chi=\one_{tC}$ for $t\in \Gamma_{/C}$. 

Now we have 
${}^{\qqq}\wtil{\Theta}(\phi, \one_{tC})_{ K, \rho} = \rr^{tC}_{ *}{}^{\qqq}\Theta(\phi)_{\rho, K}$, were $\rr^{tC}_{*}$ is induced by \eqref{dir sum}. Then 
$$\wtil{\Theta}^{tC}(\phi)_{\rho, K}:= \rr^{tC}_{*}\Theta(\phi)_{\rho, K} \quad 
 \in  H^{1}_{f}(E_{C},  M_{\rho, K})
\ot_{\Qpb} \sH_{\Qpb^{}}$$
satisfies $\underline{\qqq}(\wtil{\Theta}^{tC}(\phi)_{\rho, K}) = {}^{\qqq}\wtil{\Theta}^{tC}(\phi)_{\rho, K}$, as desired. 
\end{proof}

\subsection{The Euler system of Theta cycles}
{From now on we assume that Hyptohesis \ref{mod} holds.}

If $C\in \sC$ and  $E'$ is a finite extension of $E_{C}$, for  $z\in H^{1}_{f}(E', M_{\rho,K}^{C})$ and $f\in \sg$ we denote
$$z.f:=f_{*}z \in H^{1}_{f}(E', \rho).$$

\subsubsection{Theta cycles} \lb{sec: tc}
 For a relevant representation $\pi'\subset \sH_{\Qpb}$, denote  by $\Phi\mapsto \Phi_{\pi'}$ the Hecke-eigenprojection $\sH_{\Qpb}\to \pi'$, and by  $\lan\, , \, \ran_{\pi'^{}}\colon \pi'^{\vee}\ot\pi'\to \Qpb$ the canonical duality. We also abbreviate $\lan \vphi', \Phi_{}\ran_{\pi'} := \lan\vphi', \Phi_{\pi'}\ran_{\pi'} $ for $\vphi'\in \pi'$, $\Phi\in \sH_{\Qpb}$, and use the same names for any base-change.

For  every $\vphi\in \pi^{\vee}$, $f\in \sg$, $C\in \sC$,  and every locally constant function $\chi \colon \Gamma_{/C}\to \Qpb$, 
we define  
\beqq
   \Theta(\vphi,  \phi,  \chi)_{\rho, K}&:=\lan \vphi, \Theta(\phi,\chi)_{\rho, K }\ran_{\pi}
& \in  H^{1}_{f}(E_{C}, M_{\rho, K}),\\
 \Theta( \vphi, \phi, f, \chi)&:= \Theta( \vphi,\phi, \chi)_{\rho, K}. f& \in
  H^{1}_{f} (E_{C}, \rho),
 \eeqq
 where $K\subset \H(\A^{\infty})$ is any open compact subgroup fixing $f$ and $\phi$. 
 
As usual, in the special case $\chi=\one_{C} $ for $C\in \sC$, we will put $\Theta^{C}(-):=\Theta(-,\one_{C})$.

If $\chi\in \widehat{\Gamma}$ (viewed as an automorphic character of $\rT(\A)$), let 
 $\chi_{\G}:= \chi\circ \nu_{\G}$, $\chi_{\H}:= \chi\circ \nu_{\H}$, and denote by 
$$ \sS_{\chi} (V_{\A^{\infty}}^{r})$$ 
the space $\sS (V_{\A^{\infty}}^{r})$ with $\G(\A^{\infty})\ts \H(\A^{\infty})$-action by $\omega_{\chi}:=\omega\ot\chi_{\G}^{-1}\ot\chi_{\H}$. Let 
\beq\lb{Lm chi}
\Lambda_{\rho, \chi}:=  \left(\pi\otimes \sS_{\chi} (V_{\A^{\infty}}^{r}) \ot \sg\right)_{\G(\A_{}^{\infty}) \ts\H(\A_{}^{\infty})}, \eeq
a $\Qpb$-line. For $C\in \sC$, we also put
$$\Lm_{\rho, C}:= \bigoplus_{\chi\in\widehat{ \Gamma}_{/C}} \Lm_{\rho, \chi}.$$

\begin{lemm}[Equivariance] \lb{mod'}  Let $\chi\in \widehat{\Gamma}$.  
If  Hypothesis \ref{mod} holds, the map
\beqq \Theta(\ \cdot \ , \chi)_{}\colon \pi\otimes_{} \sS (V_{\A^{\infty}}^{r}) \ot\sg  &\to H^{1}_{f}(E, \rho(\chi))\\(\vphi_{}, \phi, f)&\mapsto \Theta(\vphi, \phi,f, \chi)
\eeqq
factors through $\Lambda_{\rho, \chi}.$
\end{lemm}

It follows from the lemma that for any $C\in \sC$, the map $\Theta^{C}$ factors through $\Lm_{\rho, C}$. 

\begin{proof}
It follows from \eqref{chi act tht} that the target is $ H^{1}_{f}(E, \rho(\chi))\subset H^{1}_{f}(E_{C(\chi)}, \rho)$. 
The equivariance for the action of $\G(\A^{\infty})$ is clear. We then need  to show that for every $\phi\in  \sS (V_{\A^{\infty}}^{r})$,  $f\in \pi$, 
  $h\in \H(\A^{\infty})$,
 we have
\beq\lb{hfchi}
\Theta(h\phi, \chi)_{\rho}. hf = \chi_{\H}^{-1}(h) \Theta(\phi , \chi)_{\rho}.f \eeq
Let  $K\subset \H(\A^{\infty})$ satisfy that $f$ and $\phi$ are invariant under $K\cap h^{-1}Kh$, and  that $\nu(K)\subset C(\chi)$. 
Then, with the notation $\phi_{| T}$ of \eqref{phiT},
we have
\beqq
Z_{T}(h\phi, \chi)_{\rho, K}. hf&= \sum_{x\in K\bks V_{\A^{\infty}}^{r}} \phi_{|T}(h^{-1}x)  Z(x,\chi)_{K} h. f& \\
&= \sum_{x\in K\bks V_{\A^{\infty}}^{r}} \phi_{|T}(h^{-1}x) \chi_{\H}^{-1}(h)  Z(h^{-1}x,\chi)_{h^{-1}Kh}.f &=  \chi_{\H}^{-1}(h)Z_{T}(\phi, \chi).f,
\eeqq
where we have used \eqref{Zchi eq} and a change of variables. This proves \eqref{hfchi}.
\end{proof}

\subsubsection{Choices of  test vectors}
Denote by:
 \begin{itemize}
 \item $\wp$ a fixed set of $p$-adic places of $F$ such that for all $v\in \wp$, $v$ splits in $E$, the Galois representation  $\rho$ is Panchishkin-ordinary (Definition \ref{ord P gal})  at each place of $E$ above $v$,  and the associated representation $\pi_{v}$ satisfies the technical condition  \eqref{cond: marcil} below;
\item $S$  a fixed set of finite places of $F$, containing all the places $v\notin \wp$ such that for some place $w\vert v$, the representation $\rho_{|G_{E_{w}}}$ is ramified or $w$ is ramified over $\Q$;
\item $\sM_{1}$  the set of split places of $E$ not in $S$;
\item $\sM$  the set of subsets of  $\sM_{1}$ (we will  identify $\sM$ and $\wp$ with a set of squarefree ideals in $\sO_{F}$);
\item $\sM[\wp]$ the set of ideals of the form $m\prod_{v\in \wp} {v}^{s_{v}}$ for $m\in \sM$ and  $s=(s_{v})\in \Z_{\geq 0}^{\wp}$. 
\end{itemize}

For $v\notin S\wp$,   let $\sV_{v}\subset V_{v}$ be a self-dual hermitian lattice in $V_{v}$, let  $K_{v}^{\circ}\subset H_{v}$ be the stabliser of $\sV_{v}^{r}$, and let $U_{v}^{\circ}\subset G_{v}$ be the stabliser of $\sum_{i=1}^{r}\sO_{E_{v}}e_{i}$. 
Fix  decompositions $\pi^{\vee}=\ot'_{v}\pi_{v}^{\vee}$, $\sg=\ot'_{v}\sg_{v}$, where the restricted tensor products are with respect to some spherical  vectors $\vphi^{\circ}_{v}\in \pi_{v}^{\vee,\circ}$,  $f_{v}^{\circ}\in K_{v}^{\circ}$ for all $v\notin S$. 

We make the following choices of  test vectors in $\pi_{v}\ot \sS(V_{v}^{r})\ot \sg_{v}$ at all finite places $v$ of $F$:
\begin{itemize}
\item
for $v\notin S \wp $,  
define
\beqq
\phi_{v}^{\circ}&:=\one_{\sV_{v}^{r}} \in \sS(V_{v}^{r})^{K_{v}^{\circ}\ts U_{v}^{\circ}}, \\
\lm_{v}^{\circ}&:= \vphi_{v}^{\circ} \ot \phi_{v}^{\circ}\ot f_{v}^{\circ}.\eeqq
\item for $v\in S$, we let $\lm_{v}=\vphi_{v}\ot\phi_{v}\ot f_{v}\in \pi_{v}\ot \sS(V_{v}^{r}) \ot \sg_{v}$ be any element whose image in $\Lm_{\rho, v}$ is nonzero;
\item
for  $v\in \sM_{1}$, we 
will define another Schwartz function 
$$\phi_{v}^{\bullet}:= \eqref{phiv def} \in \sS(V_{v}^{r})^{K_{v}^{\circ}\ts U_{v}^{\circ}}$$
below (where the subscripts $v$ will be  omitted from the notation), and we put
$$\lm_{v}^{\bullet}:= \vphi_{v}^{\circ} \ot \phi_{v}^{\bullet}\ot f_{v}^{\circ};$$
\item for $v\in \wp$, we will define vectors $\vphi^{\ra}\in \pi^{\vee}$, $f^{\ra}\in \sg$ and  a sequence of Schwartz functions
$$\phi_{v}^{(s)}\in \sS(V_{v}^{r}), \qquad s\geq 0,$$
in Definition \ref{test p} below (where the subscripts $v$ will be  omitted from the notation). We put 
$$\lm_{v}^{(s)}:= \vphi^{\ra}\ot \phi_{v}^{(s)}\ot f^{\ra}.$$
\end{itemize}
For $m=\prod_{v|m, v\in \sM_{1}} v\prod_{v\in \wp} v^{s_{v}}\in \sM[\wp] $, we put
$$  \lm^{(m)} := (\ot_{v\notin S m\wp} \lm_{v}^{\circ})\ot ( \ot_{v|m, v\in \sM_{1}} \lm_{v}^{\bullet}) \ot (\ot_{v\in \wp} \lm_{v}^{(s_{v})})\ot \ot_{v\in S}\lm_{v} \qquad \in \pi^{\vee}\ot \sS(V_{\A^{\infty}}^{r})\ot \sg; $$
we also put $\lm^{(0)}:=\lm^{(1)}$, and define $\vphi^{(m)}\in \pi^{\vee}$, $f_{v}^{(m)}\in \sg$ in the obvious way so that $\lm^{(m)}=\vphi^{(m)}\ot \phi^{(m)}\ot f^{(m)}$. Note that $\vphi^{(m)}$ and $f^{(m)}$ are in fact independent of $m\in \sM[\wp]\cup\{0\}.$

\subsubsection{The Euler system} \lb{sec: ES}
For $m\in \sM[\wp]$,  we set 
 $$C(m):=  (1+m\widehat{\sO}_{E})\cap \rT(\A^{\infty}), \qquad E[m]:=E_{C(m)}.$$
For $m=0$, we put  $C(0)=  \rT(\A^{\infty})$, $E[0]:=E$. 

\begin{defi}\lb{def th c}
The \emph{Euler system of Theta cycles} is the system of classes $(\Theta_{m})_{m\in \sM[\wp]\cup\{0\}}$ defined by 
\beqq
\Theta_{m}:= 
 \Theta^{C(m)}( \lm^{(m)})\qquad
\in H^{1}_{f}(E[m], \rho) .\eeqq
\end{defi}

The following theorem says precisely that $(\Theta_{m})_{m}$ is an Euler system in the sense of Jetchev-\nek-Skinner.
For a place $w$ of $E$ at which $\rho$ is unramified, let $\r{Fr}_{w}\in G_{E_{w}}$ be a geometric Frobenius at $w$, and let $P_{w}(t):= \det(1- t \r{Fr}_{w}|\rho^{*}(1))$.  

\begin{theo}  \lb{norm rel}
The system of classes $(\Theta_{m})_{m\in \sM[\wp]\cup \{0\}}$ of  Definition \ref{def th c}
 satisfies $\Tr_{E[1]/E}\Theta_{1}=\Theta_{0}$ and the following conditions.
\begin{enumerate}
\item\lb{thm int}
 \emph{Integrality.}  There exists a $G_{E}$-stable $\Zpb$-lattice $\rho_{0}\subset \rho$ such that for every $m\in \sM[\wp]\cup \{0\}$, 
$$\Theta_{m} \in H^{1}_{f}(E[m], \rho_{0}).$$
\item\lb{NRh}
  \emph{Horizontal norm relations.} For every $m\in \sM$ and every $v\in  \sM_{1}$ not dividing $m$,
\beqq
\Tr_{E[mv]/E[m]} \Theta_{mv} = P_{w}(\Fr_{w})\,  \Theta_{m},\eeqq
where $w$ is any one of the places of $E$ above $v$.
\item\lb{NRv}
 \emph{Vertical norm relations.} For every $m\in \sM[\wp]$ and every $v\in \wp$, 
\beqq
\Tr_{E[mv]/E[m]} \Theta_{mv} =   \Theta_{m}.\eeqq
\end{enumerate}
\end{theo}

\begin{rema} 
By construction, the $\Qpb$-vector space $\Lambda_{\rho, S}=\ot_{v}\Lm_{\rho,v}$ is $1$-dimensional; thus the `base class' $\Theta_{0}$, which only depends on   the image of $\lm_{S}=\ot_{v}\lm_{v}$ in $\Lm_{\rho, S}$, is independent of choices   up to a scalar (after the initial choice of the descent $\pi$).  The following proposition verifies the resulting necessary condition for the  nonvanishing of (the base class of) our Euler system.
\end{rema}

We say that $\rho$ is \emph{exceptional} at a place $v\in \wp$ if for some (equivalently,\footnote{This follows from multiplicativity and functional equation of $\gamma$-factors, the selfduality of $\rho$, and the fact that, by weight considerations, for every $w$ the factor $\gamma(\mathrm{WD}_{\iota}(\rho_{w}),\psi_{E, w}, s)$ has neither a zero nor a pole at $s=0$.}
 every) place $w\vert v$ of $E$ and embedding $\iota\colon \Qpb\into \bC$, the Deligne--Langlands $\gamma$-factor 
$$\gamma({\mathrm{WD}}_{\iota}(\rho_{w}^{+}),\psi_{E,w}, s)$$
of the complex Weil--Deligne representation attached to $\rho_{w}^{+}$ by \cite{Fon94} does not have a pole at $s=0$. A consideration of weighs shows that if  $\rho$ is crystalline at all $w\vert v$, then it is not exceptional. 
\begin{prop}\lb{NV} The image of  $\lm^{(0)}$ in $\Lm_{\rho}$ is nonzero if and only if $\rho$ is not exceptional at any place $v\in \wp$.
\end{prop}
This is clearly a local statement, which will be proved in \S~\ref{loc nv}. 

\subsection{Reduction of the Euler-system properties to local statements} 
We reduce Theorem \ref{norm rel} to several local results, to be proved in the next section; for  clarity, these results are marked with a `$\to$'.  

\subsubsection{Integral structures} Let $K\subset \H(\A^{\infty})$, $U\subset \G(\A^{\infty})$  be   compact open subgroups fixing $\vphi^{(m)}$ and $f^{(m)}$. We consider the following integral structures on our representations.
\begin{itemize}
\item  We let $\rho_{\Zpb}$ be a $\Zpb$-lattice in $\rho$, stable under $G_{E}$ (this may require a choice that we now fix);
\item Let $M_{\rho,\Zpb, K}\subset M_{\rho, K}$ be a $\Zpb$-lattice such that for each $C$, the image of $\Ch_{r-1}(X_{K, E_{C}})_{\Zpb}$ of the cycle class map $(-)_{\rho}$ from Lemma \ref{p rho} is contained in  $H^{1}(E_{C}, M_{\rho, \Zpb, K})$. (As explained in 
  \cite[\S~II.1.10]{nek-heeg},  we may  take $M_{\rho, \Zpb, K}= p^{-a} e_{\rho} H^{2r-1}_{\text{\'et}}(X_{K, \ol{E}}, \Zpb(r))$, where $p^{a}$ is the order of the torsion subgroup of $H^{2r-1}_{\text{\'et}}(X_{K, \ol{E}}, \Zpb(r))$.)
  Let $M_{\sg, \Zpb, K}
:= M_{\sg, K}\cap M_{\rho, \Zpb, K}$. 
We define
$$\sg_{\Zpb}^{K}:= \Hom_{\Qpb[K\bks \H_{V}(\A^{\infty}) / K]}( M_{\sg, \Zpb, K}, \rho_{\Zpb}),$$
a $\Zpb$-lattice.
\item Let $\sH_{\Zpb}$ be the preimage of $\ul{\r{SF}}_{\Zpb}$ under the map $\ul{\qqq}$ of \eqref{ul q}. We define 
$$\pi^{U,\vee}_{\Zpb}:=\pi^{U,\vee}\cap \Hom(\sH_{\Zpb}^{U}, \Zpb),$$ where $\pi^{U,\vee}$ is viewed as a subspace of $\Hom(\sH_{\Qpb}, \Qpb)$ via the composition of the natural duality and the projection $\sH_{\Qpb}\to \pi$. 
\item For $x\in V_{\A^{\infty}}^{r}$ and $C\in \sC$, let $V(x)$ and $\H(x)(\A^{\infty})$ be as in \S~\ref{sec: sp cy}, and let $K_{x}:= K\cap \H(x)$, $K_{x}^{C}:= K_{x}\cap \nu_{\H(x)}^{-1}(C)$. For each $C\in \sC$, we define
 $$ \sS(V_{\A^{\infty}}^{r}, \ol{\mathbf{Z}}_{p,C})^{K} \subset  \sS(V_{\A^{\infty}}^{r}, \Qpb)^{K}$$ 
 to be the $\Zpb$-module of functions satisfying that for every $x\in \mathrm{Spt}(\phi)$, 
$$\vol(K) \cd \phi(x) \cd [K_{x}:K_{x}^{C}]^{-1} \in \Zpb.$$
Similar integrality properties for Schwartz functions  are considered by Shah in \cite[\S 3.5]{Shah}. 
\end{itemize}

\begin{rema}\lb{lattice}
The $\Zpb$-submodule  $\sH_{\Zpb}^{U}$ is a $\Zpb$-lattice in the subspace $\sH_{\Qpb}^{U, \circ}\subset \sH_{\Qpb}^{U}$ consisting of forms with (uniformly) bounded $q$-expansions. (We conjecture that $\sH_{\Qpb}^{U, \circ}= \sH_{\Qpb}^{U}$; at least when $U$ is hyperspecial at $p$-adic places, this should be provable by considering $q$-expansion maps on integral models of PEL Shimura varieties related to $\G$, cf. \cite[Remark 5.2.14]{Lan12}.)  
This implies that  $\pi_{\Zpb}^{U,\vee}\subset \pi^{U, \vee}$ contains a $\Zpb$-lattice.
\end{rema}

\begin{lemm}\lb{suff int}
 For every $C\in \sC$, we have 
\beqq
\vphi\in \pi^{U, \vee}_{\Zpb},  \quad 
 \phi\in \sS(V_{\A^{\infty}}^{r}, \ol{\mathbf{Z}}_{p,C}))^{U\ts K}, \quad f\in \sg^{K}_{\Zpb}
 \qquad \Longrightarrow \qquad
\Theta^{C}( \vphi, \phi, f)\in H^{1}_{f}(E_{C}, \rho_{\Zpb}).\eeqq
\end{lemm}
\begin{proof}
By the definitions, it suffices to prove that 
$$\Theta^{C}( \phi)_{K}\in  \Ch_{r-1}(X_{K, E_{C}})_{\Zpb}\ot_{\Zpb} \sH_{\Zpb}^{U},$$
that is, that for all $x\in {\r{Spt}}(\phi)$ with $T(x)\in \Herm_{r}(F)$ and for all $t\in T(\A^{\infty})$, we have
$$ \vol(K)\cd \phi(x) \cd [K_{x}: K_{x}^{C}]^{-1}\cd 
\left(c_{1}(\sL_{K}^{\vee})^{\dim V(x)-r}_{|X^{tC}_{K}}
\frown[j_{x}(X_{\rH(x), K_{x}}^{tC})]  \right) \in  \Ch_{r-1}(X_{ K, E_{C}})_{\Zpb}.$$
This is immediate from the definition of $\sS(V_{\A^{\infty}}^{r}, \ol{\mathbf{Z}}_{p,C}))^{U\ts K}$.
\end{proof}

\subsubsection{Integrality} \lb{inttt}   We reduce the integrality of the Euler system to a result on the local integrality of our Schwartz functions.

For $x_{v}\in V_v^{r}$ , let $V(x_{v})=\Span(x_{v,1}, \ldots, x_{v,r})^{\perp}$, let $H(x_{v})= U(V(x_{v}))$. For open commpact subgroups $U_{v}\subset G_{v}$,  $K_{v}\subset H_{v}$ and $C_{v}\subset T_{v}$, let  $K_{x_{v}}:= K_{v}\cap H(x_{v})$, $K_{x_{v}}^{C_{v}}:= K_{x_{v}}\cap \nu_{H(x_{v})}^{-1}(C_{v})$. 
We define
 $$ \sS(V_v^{r}), \ol{\mathbf{Z}}_{p,C_{v}})^{U_{v}\ts K_{v}} \subset  \sS(V_v^{r})^{U_{v}\ts K_{v}}$$ 
to  be the  $\Zpb$-module of functions $\phi_{v} $ satisfying that for every $x_{v}\in \mathrm{Spt}(\phi_{v})$, 
$$\vol(K_{v}) \cd \phi_{v}(x_{v}) \cd [K_{x_{v}}:K_{x_{v}}^{C_{v}}]^{-1} \in \Zpb.$$

For a finite place $v$ of $F$, we denote by $\vpi_{u}$ a fixed uniformiser. If $v$ is a split finite place of $F$, we  denote  $C(\vpi_{v}^{e}) :=1+\vpi_{v}^{e} \sO_{E,v} \cap T_{v}$.  
 The following will be proved as Propositions  \ref{inttt2}   and \ref{intt p} below.
\begin{enonce}{Proposition${}^\to$} \lb{int}   We have:
\begin{enumerate}
\item \lb{int M} for $v\in \sM_{1}$, 
$$\phi_{v}^{\bullet}\in  \sS(V_v^{r}), \ol{\mathbf{Z}}_{p,C(\vpi_{v})})^{U_{v}^{\circ} \ts K_{v}^{\circ}};$$
\item \lb{int p} for $v\in \wp$, there exist open compact subgroups $U_{v}\subset G_{v}$, $K_{v}\subset H_{v}$  fixing respectively $\vphi^{\ra}_{v}$ and $f^{\ra}_{v}$, such that for every $s\geq0$, 
$$\phi_{v}^{(s)}\in  \sS(V_v^{r}, \ol{\mathbf{Z}}_{p,C(\vpi_{v}^{s}})^{U_{v}\ts K_{v}}.$$
\end{enumerate}
\end{enonce}

\begin{proof}[Proof of Theorem \ref{norm rel}.1, assuming Proposition \ref{int}]
Proposition \ref{int} implies that there are compact open subgroups $U\subset \G(\A^{\infty})$ and $K\subset \H(\A^{\infty})$,  and an integer $e_{2}$ such that for every $m\in \sM[\wp]$, 
$$\phi^{(m)}\in p^{-e_{2}}\sS(V_{\A^{\infty}}^{r}, \ol{\mathbf{Z}}_{p,C(m)})^{K\ts U}.$$
Now  Lemma \ref{suff int} shows that if  $e_{1}, e_{3}\in \Z$ are such that $\vphi^{(m)}\in p^{-e_{1}} \pi^{\vee, U}_{\Zpb}$, $f^{(m)}\in p^{-e_{3}}\sg^{K}_{\Zpb}$  (recall that $\vphi^{(m)}$ and $f^{(m)}$ are independent of $m$, so this is possible by Remark \ref{lattice}), then for the $\Zpb$-lattice
$$\rho_{0}:=p^{-e_{1}-e_{2}-e_{3}}\rho_{\Zpb}\subset \rho$$
we have $\Theta_{m}\in H^{1}_{f}(E[m], \rho_{0})$, as desired.
\end{proof}

\subsubsection{Horizontal norm relations}

We first reduce the norm relation of Theorem \ref{norm rel}.\ref{NRh} to the following proposition.

\begin{enonce}{Proposition${}^\to$} \lb{norm rel'} 
For all $m\in \sM$ and  $v\in \sM_{1}$ with $v\nmid m$,  all characters  $\chi \colon {\Gamma_{/C(m)}}\to \Qpb^{\ts}$, 
and for every $\lm_{v}\in  \pi^{\vee,v}\ot \sS(V_{\A^{v\infty}}^{r})\ot \sg^{v}$, we have
\beq\lb{eq: nr'}
\Theta(\lm^{v}\lm_{v}^{\bullet}, \chi)_{}= P_{w}({\mathrm{Fr}}_{w})
\Theta(\lm^{v}\lm_{v}^{\circ},\chi)_{}\eeq
in $H^{1}_{f}(E[m], \rho)\ot_{\Qpb} \sH_{\Qpb}$.
\end{enonce}

\begin{proof}[Proof of Theorem \ref{norm rel}.\ref{NRh}, assuming Proposition \ref{norm rel'}]
The identity  \eqref{eq: nr'} remains valid for any  function $\chi$ on $ {\Gamma}_{/C(m)}$.  
Then it suffices to apply it to $\chi=\one_{C(m)}$ and $\lm^{v}=\lm^{(m),v}$. 
\end{proof}

We can further reduce Proposition \ref{norm rel'}  to the following abstract local analogue, 
to be proved in \S~\ref{sec loc nrh}.
Analogously to \eqref{Lm chi}, denote by $\sS_{\chi_{v}}(V_{v}^{r})$ the space $\sS(V_{v}^{r})$ equipped with the $(G_{v}\ts H_{v})$-action by $\omega_{v, \chi_{v}}:=\omega_{v}\ot\chi_{\G, v}^{-1}\ot\chi_{\H, v}$, and let 
\beq\lb{Lm chi v}
\Lm_{\rho_{v}, \chi_{v}}^{}:= 
\left( \pi_{v}^{\vee}\ot_{}\sS_{\chi_{v}}(V_{v}^{r}) \ot \sg_{v}\right)_{G_{v}\ts H_{v}}. \eeq
\begin{prop} \lb{norm rel''} For every  $v\in \sM_{1}$ and every unramified character $\chi$ of $E_{v}^{\ts}$, we have
\beq\lb{tht id}
[\lm_{v}^{\bullet}]
= L(\rho^{*}(1)_{w}, \chi^{-1}_{w}, 0)^{-1}\cd 
[\lm_{v}^{\circ}]
\eeq
in 
$\Lm_{\rho, \chi, v}$.
\end{prop}
\begin{proof}[Proof of Proposition  \ref{norm rel'}, assuming Proposition \ref{norm rel''}]
Characters $\chi\colon \Gamma_{/C(m)}\to \Qpb^{\ts}$ are 
unramified at $v\nmid m$, and  by   \eqref{chi act tht},
$\Fr_{w}$ acts by $\chi^{-1}(\mathrm{Fr}_{w})$ on  the generating series $\Theta(\phi, \chi)_{\rho}$. Thus the Galois element $P_{w}(\Fr_{w})$ acts by  the scalar $P_{w}(\chi^{-1}(\mathrm{Fr}_{w}))=L(\rho^{*}(1)_{w},\chi_{w}^{-1},0)^{-1}$, and  the desired identity  \eqref{eq: nr'} simplifies to 
$$\Theta(\lm^{v}\lm_{v}^{\bullet}, \chi)
=L(\rho^{*}(1)_{w}, \chi_{w}^{-1},0)^{-1}\cd 
\Theta(\lm^{v} \lm_{v}^{\circ},\chi).$$
This identity is implied by Proposition \ref{norm rel''} since, by Lemma \ref{mod'}, the map $\lm_{v}\mapsto\Theta(\lm^{v}\lm_{v}, \chi)$ factors through $\Lm_{\rho, \chi,v}$.
\end{proof}

\subsubsection{Vertical norm relations}
We reduce Theorem \ref{norm rel}.\ref{NRv} to the following.
\begin{enonce}{Proposition${}^\to$}\lb{prop NRv}
Let $v\in \wp$. For every compact open subgroup $C_{v}\subset E_{v}^{\ts}$, the image of $\lm_{v}^{(s)}$ in $\Lm_{\rho_{v}}^{C_{v}}=\bigoplus_{\chi_{v}\in \widehat{E_{v}^{\ts}/C}}\Lm_{\rho_{v}, \chi_{v},v}$ is independent of $s\geq 0$. 
\end{enonce}

\begin{proof}[Proof of Theorem \ref{norm rel}.\ref{NRv}, assuming Proposition \ref{prop NRv}] 
Let $s=v(m)$. We have
$$\Tr_{E[mv]/E}\Theta_{mv}=\Tr_{E[mv]/E[m]}\Theta^{C(mv)}(\lm^{(m), v}\lm_{v}^{(s+1)})=\Theta^{C(m)}(\lm^{(m), v}\lm_{v}^{(s)})=\Theta_{m}.$$
\end{proof}

\section{Local study}\lb{sec 3}

The goal of this section is to prove Propositions \ref{int}, \ref{norm rel''}, and \ref{prop NRv}, after giving the definition of the test vectors.

\Subsection{Preliminaries}

\subsubsection{Notation}
Let $v $ be a finite place of $F$ split in $E$. We work in a local setting over $F_{v}$ and drop all subscripts $v$ (thus writing $F$, $E$, $V$, $W$, \ldots\  for $F_{v}$, $E_{v}$, $V_{v}$, $W_{v}$, \ldots).   We denote by $d\in \sO_{F}$ a generator of the different ideal of $F$. 

 We fix an ordering of the two primes $w_{1}, w_{2}$ above $v$ and write $E=E_{w_{1}}\times E_{w_{2}}=F\ts F$; we put $\ri:= (1, -1)\in E$.  We fix a uniformiser $\vpi$ of $F$, and we denote by ${k}$ the residue filed of $F$,  by $q$ its cardinality.

Write $V=V_{v}= V_{1}\oplus V_{2}$ where  the isotropic subspaces $V_{i}:=V\ot_{F} E_{w_{i}}$. Write $W=E^{2r}=W^{+}\oplus W^{-}$  where $W^{+}=\Span(e_{1}, \ldots, e_{r})$, $W^{-}=\Span(e_{r+1}, \ldots, e_{2r})$. For $?=\emptyset, +, -$, write $W^{?}=W^{?}_{1}\oplus W^{?}_{2}$.   Let $\sW:= \bigoplus_{i=1}^{n} \sO_{E} e_{i}$, and if `$?$' is any decoration, let $\sW_{?}=W_{?}\cap \sW$.

We denote
$$H=\H(F)\cong {\mathrm{Aut}}_{F}(V_{1}), \qquad G=\rG(F)\cong {\mathrm{Aut}}_{F}(W_{1}),\qquad T= \rT(F)\cong F^{\ts},$$ 
where $H$ acts on $V_{1}$ on the left, and $G$ acts on $W_{1}$ on the right. 

Fix an isometry between $V$ and $W^{\vee}$, where the latter is endowed with the hermitian form dual to the  form $(y, y')_{W}:=\ri  \lan y, y'\ran_{W}$ on $W$. 
The chosen isometry between $V$ and (the hermitian space attached to) $W^{\vee}$ induces isomorphisms\footnote{Note that $G$ acts on the left on $W_{1}^{\vee}$.}
\beq\lb{id g h}
H=\mathrm{Aut}_{F}(V_{1})  \to\mathrm{Isom}_{F} ( W_{1}^{\vee}, V_{1})  \stackrel{}{\leftarrow} \mathrm{Aut}_{F}(W_{1}^{\vee})=G.\eeq

By the isomorphism $W_{i}\cong V_{i}^{\vee} $ (for $i=1,2$), we have $\sO_{F}$-lattices $\sV_{i}=\sW_{i}^{\vee}$, and   direct summands $V_{i}^{\pm} = (W_{i}^{\pm})^{\vee}\subset V_{i}$; we let $\sV_{i}^{\pm}:= \sV_{i}\cap V_{i}^{\pm}$. We will often write
$$x_{i}^{\pm} = {x_{i, +}^{\pm}\choose  x_{i, -}^{\pm} } \in V_{i}\ot W_{i}^{\pm}
= \Hom(V_{i}^{\pm}, V_{i}) 
=\Hom(V_{i}^{\pm}, V_{i}^{+})
\oplus 
\Hom(V_{i}^{\pm}, V_{i}^{-}).
$$ 

\subsubsection{Subgroups of $G$ and $H$}
We denote still by  $\rH$  the algebraic group over $\sO_{F}$ with $\H(R) =\mathrm{Aut}_{R}(\sV_{1}\ot_{\sO_{F}}R)$ for any $\sO_{F}$-algebra $R$; we similarly extend $\G$ to a group over $\sO_{F}$ (isomorphic to $\rH$). 
 We write $P=\rP(F)\subset G=\GL(W_{1})=\GL_{n}(F)$ for the Siegel parabolic, with Levi $M\cong \GL(W_{1}^{+})\ts \GL(W_{1}^{-}) =: G^{+}\ts G^{-}$ and unipotent radical $N=\rN(F)$. 
 
We put
$$K^{\circ}:=\G(\sO_{F}),$$ 
and define (deeper) 
 pro-$p$ parahoric subgroups of  $\G(\sO_{F})$  of level $s\geq 1$ by
 $$
 \quad  I_{s} = \rG(\sO_{F})\ts_{\rG(\sO_{F}/\vpi^{s}\sO_{F})} \rN(\sO_{F}/\vpi^{s}\sO_{F}).$$
 Via the identification \eqref{id g h}, we may also view the above as subgroups of $H$; when the context is ambiguous (and the distinction is needed), we will add a superscript `$G$' or `$H$' to the notation for those subgroups in order to distinguish the ambient group in question.

Finally, we will need the following definition. 
\begin{defi}\lb{gal level}
Let $H'$ be a general linear group over $F$ and  let $s\geq 0$ be an integer. We say that a subgroup $K\subset H'$ has \emph{Galois-level at least $s$} if $\det( K) \subset 1+\vpi^{s}\sO_{F} $.
\end{defi}

For the rest of this section, unless noted otherwise: all tensor products of finite-dimensional $F$-vector spaces are taken over $F$; Schwartz spaces consist of $\Qpb$-valued functions; all tensor products of $\Qpb$-vector spaces are taken over $\Qpb$.

\subsubsection{Weil action and linear action on Schwartz spaces} For any smooth character of $F^{\ts}$, denote by 
$$ \sS_{\chi}'(V_{1}\ot W_{1}) $$
the Schwartz space of $V_{1}\ot W_{1}$
endowed with the action of  $G\ts H$ given by 
\beq \lb{lin act}
(g, h).\phi'(y)= \chi^{-1}(\det g)\chi(\det h)|\det h|^{-r} |\det g|^{r} \phi'( h^{-1}y  g).\eeq
We still denote by $\sS_{\chi}(V\ot_{E} W^{+})$ the Schwartz space of 
$$V\ot_{E} W^{+} \cong 
  (V_{1}\ot W_{1}^{+})\oplus (V_{2}\ot W_{2}^{+}) $$
endowed with the twisted  Weil action $\omega_{\chi}=\omega \ot \chi_{G}^{-1}\ot \chi_{H}$. (When we are interested in the Schwartz space only and not the specific action, we will omit the subscript $\chi$.)

Define a linear map
\beq
\sF\colon \sS_{\chi}'(V_{1}\ot W_{1}) =\ & \sS_{\chi}'((V_{1}\ot W_{1}^{+}) \oplus (V_{1}\ot W_{1}^{-}) ) \\
\to \ &
\sS_{\chi} ( (V_{1}\ot W_{1}^{+}) \oplus (V_{2}\ot W_{2}^{+}))=\sS_{\chi}(V\ot_{E} W^{+})
\eeq
 by the partial Fourier transform
$$\sF\phi'(x_{1}^{+}, x_{2}^{+}) = \int_{V_{1}\ot W_{1}^{-}}  \phi'(x_{1}^{+}, x_{1}^{-}) \psi(\lan x_{1}^{-}, x_{2}^{+}\ran) \, \rd x_{1}^{-}
$$
where $\lan\, , \, \ran$ is the natural duality between $V_{1}\ot W_{1}^{-}$ and $V_{2}\ot W_{2}^{+}$ given by the restriction of  $( \ ,\ )_{V}\ot (\  ,\  )_{W}$, and $\rd x_{1}^{-}$ is the self-dual Haar measure, which assigns volume $|d|^{r^{2}}$ to $\sV_{1}\ot_{\OO_{F}} \sW_{1}^{+}$. 

\begin{lemm} \lb{intertw}The map $\sF \colon \sS_{\chi}'(V_{1}\ot W_{1}) \to \sS_{\chi}(V\ot_{E} W^{+}) $ is an isomorphism of $(G\ts H)$-modules.
\end{lemm}
\begin{proof} 
It is easy to verify by explicit computation that $\sF$ is $(H\ts G)$-equivariant (for a brief discussion in a more general context, see  \cite[\S~ 2.9]{GQT}). It is also clear that the dual partial Fourier transform gives an explicit inverse. 
\end{proof}

\subsubsection{Godement--Jacquet zeta integrals as models for  the Howe correspondence}
We fix an embedding $\Qpb\into \bC$ extending $\iota^{\cyc}$, via which we may base-change all $L$-values, functionals, and representations -- without changing the notation. 

 By the uniqueness of  the Howe correspondent $\sg$ in \S~\ref{p aut rep} together with Lemma \ref{GJ gen} below, we have $\sg\cong \pi$ under the isomorphism \eqref{id g h}.
 Fixing a nontrivial $H$-equivariant  pairing $(\ , \ )_{\pi}\colon \pi^{\vee} \ot \pi \to \bC$, 
we have the Godement--Jacquet zeta integral 
\beq 
\zeta(\cd, \chi) \colon  \pi^{\vee}\ot \sS'_{\chi}(V_{1}\ot W_{1})  \ot \sg & \to \bC\\
 ( \vphi, \phi', f ) & \mapsto
  \int_{G} (g\vphi, f)_{\pi} \cd \phi'(
g)\chi^{-1}(\det g)  |\det g|^{r} \, \rd g
 \eeq
where $\rd g $ is the Haar measure on $G$ assigning volume~$|d|^{n^{2}/2}$ to $K^{\circ}:=\mathrm{Aut}_{\sO_{F}}(\sV_{1})\subset G\cong H$. 
Let 
\beq \lb{def tht}
\theta(\cd, \chi) := \zeta(\cd, \chi) \circ \sF^{-1}\in \Lm_{\rho, \chi}^{\vee}\eeq
(where $\Lm_{\rho, \chi}=\eqref{Lm chi v}$).
\begin{lemm} \lb{GJ gen}
The functional $\theta(\cd, \chi)$ is a generator of  $\Lm_{\rho, \chi}^{\vee}$. 
\end{lemm}
\begin{proof} By Lemma \ref{intertw}, this is equivalent to the assertion that $\zeta(\cd, \chi)$ is  a nonzero element of $\Hom_{G\ts H}(\pi^{\vee}\ot \sS_{\chi}'(V_{1}\ot W_{1}) \ot\sg, \bC)$. The belonging is easily verified. For the nonvanishing, it suffices to apply $\zeta(\cd, \chi)$ to a triple $(\vphi, \phi', f)$ such that $(\vphi, f)_{\pi}\neq 0$ and $\phi'$ has  small support near the identity.
\end{proof}

\subsection{Test vectors, norm relations and integrality at places in $\sM_{1}$}\lb{sec: horiz}
Suppose that $v\in \sM_{1}$.

\subsubsection{Definition of the Schwartz function $\phi^{\bullet}$}
If $P$ is a logical proposition,  denote by $$\one[P]\in \{1=\text{true}, 0=\text{false}\}$$ its truth value (thus for the characteristic function of a set $A$ we have $\one_{A}(x)=\one[x\in A]$). 
 Let
\beqq 
\phi'^{\bullet}
:=  \one_{\mathrm{Aut}_{\sO_{F}}(\sV_{1})} = \sum_{h\in K^{\circ}/I_{1}} h\phi_{0}'^{(v)} \qquad \in \sS'(V_{1}\ot W_{1}), 
\eeqq
where 
$$\phi_{0}'^{\bullet}\left(x_{1} \right):= 
\one \left[x_{1}=  \twomat  {x_{1, +}^{+} }  {x_{1, +}^{-}}  {x_{1, -}^{+} } {x_{1, -}^{-} }
\in \twomat 
{\id_{V_{1}^{+}} +\vpi \End(\sV_{1}^{+})}
{\Hom(\sV_{1}^{-},  \sV_{1}^{+})}
{\vpi\Hom(\sV_{1}^{+}, \sV_{1}^{-})}
{\id_{V_{1}^{-}}+\vpi \End(\sV_{1}^{-})}
\right].
$$
The  stabiliser of $\phi_{0}'^{\bullet}$ under the action 
 \eqref{lin act} (for $\chi=\one$) of $K^{\circ} \subset H={\mathrm{Aut}}_{F}(V_{1})$  is $I_{1}$.

W define 
\beq \lb{phiv def}
\phi^{\bullet}:= \sF \phi'^{\bullet} = \sum_{h\in K^{\circ}/I_{1}} h\phi_{0}^{\bullet}, \qquad \phi_{0}^{\bullet}:=\sF\phi_{0}'^{\bullet} .\eeq
The choice of $\phi^{\bullet}$ is motivated by Proposition \ref{GJ calc} below.

\subsubsection{Local horizontal norm relations} \lb{sec loc nrh}
We  prove Proposition \ref{norm rel''} by a computation of zeta integrals.

\begin{prop}\lb{GJ calc} 
 For $?\in\{\circ, \bullet\}$, let $\lm^{?}:= \vphi^{\circ}\ot \phi^{?} \ot f^{\circ} \in \pi \ot\sS(V^{r})\ot \sg$.
 For every unramified character $\chi$ of $F^{\ts}$, we have
$$
[\lm^{\bullet}]= L(1/2, \pi(\chi)^{\vee})^{-1}\cd [\lm^{\circ}] $$
in $ \Lm_{\rho, \chi}$.
\end{prop}
\begin{proof}  
 By \cite[Lemma 6.10]{GJ} and, respectively, the definition,
 we have
 \beqq
 \zeta (\vphi^{\circ}, \phi'^{\circ}, f^{\circ},\chi)&=  L(1/2, \pi(\chi)^{\vee}) \cd  (  \vphi^{\circ}, f^{\circ})_{\pi}, \\
\zeta (\vphi^{\circ}, \phi'^{\bullet}, f^{\circ}, \chi)&= ( \vphi, f)_{\pi}.
\eeqq 
By Lemma \ref{GJ gen}, this implies the desired result.
\end{proof}
\begin{proof}[Proof of   Proposition \ref{norm rel''}]
It is equivalent to  Proposition \ref{GJ calc}, once noted that $L(\rho^{*}(1)_{w}, \chi_{w}^{-1}, s)=L(s+1/2, \pi(\chi)_{v}^{\vee})$.\end{proof}

\subsubsection{A decomposition}
We begin a study of the function $\phi^{\bullet}$, with the final goal of establishing its   integrality properties. 

We have
\beqq
\sF\phi_{0}'^{\bullet} (x_{1}^{+}, x_{2}^{+})
&=
\one \left[\twomat  {x_{1, +}^{+} }  {x_{2, +}^{+}}  {x_{1, -}^{+} } {x_{2, -}^{+} }
\in \twomat 
{\id_{V_{1}^{+}} +\vpi \End(\sV_{1}^{+})}
{\vpi^{-1} \End(\sV_{2}^{+})}
{\vpi\Hom(\sV_{1}^{+}, \sV_{1}^{-})}
{ \Hom(\sV_{2}^{+}, \sV_{2}^{-})}
\right]
\cd
{\psi}(\Tr(x_{2, +}^{+}))\\
&=
 \sum_{y\in Y}
\phi_{0, y}^{\bullet}(x_{1}^{+}, x_{2}^{+}), 
\eeqq
where $Y:= \vpi^{-1}  \End(\sV_{2}^{+})/  \End(\sV_{2}^{+})$ and 
\beqq
\phi_{0, y}^{\bullet}(x_{1}^{+}, x_{2}^{+}):= 
{\psi}(\Tr (y))\cd
\one \left[\twomat  {x_{1, +}^{+} }  {x_{2, +}^{+}}  {x_{1, -}^{+} } {x_{2, -}^{+} }
\in \twomat 
{\id_{V_{1}^{+}} +\vpi \End(\sV_{1}^{+})}
{y+ \End(\sV_{2}^{+})}
{\vpi\Hom(\sV_{1}^{+}, \sV_{1}^{-})}
{ \Hom(\sV_{2}^{+}, \sV_{2}^{-})}
\right].
\eeqq

\subsubsection{Rationality} It is clear that the functions $\phi_{0, y}^{\bullet}$ take values in $\Q(\mu_{\ell^{\infty}})$, where $\ell$ is the rational prime underlying $v$. The following lemma is not strictly necessary, but we include it for completeness.
\begin{lemm}  \lb{rat}
For $u\in \Z_{\ell}^{\ts}$, let $\sg_{u}\in \Gal(\Q(\mu_{\ell^{\infty}}/\Q))$ be its image under the reciprocity map. 
\begin{enumerate}
\item
For all $u\in \Z_{p}^{\ts}$, we have $\phi_{0,y}^{\bullet}(x)^{\sg_{u}} = \phi^{\bullet}_{0, u^{-1}y}\left( \smalltwomat {1_{r}}{}{}{u 1_{r}} x\right)$. 
\item The function $\phi^{\bullet}$ takes values in $\Z$. 
\end{enumerate}
\end{lemm}
\begin{proof}
Part 1  follows  from the relation  $\psi(t)^{\sg_{u}}=\psi(u^{-1}t)$ (for all $t\in F$) and the definitions, using the description of the action of  $ \smalltwomat {1_{r}}{}{}{u 1_r}\in \GL_{n}(F)\cong H$ on $V\ot_{E} W^{+}$ given in \eqref{h act}, \eqref{h un} below. From part 1 and  the definitions, since $\smalltwomat {1_{r}}{}{}{u 1_{r}}\in K^{\circ}$ it follows that $\Gal(\Q(\mu_{\ell}^{\infty})/\Q))$ fixes $\phi^{\bullet}$; it is also clear that the values of $\phi$  are integers.
\end{proof}

\subsubsection{Levels and integrality}
We now consider  $V$ as endowed with the dual basis $b_{1},\ldots, b_{n}$ to the standard basis  of $W=E^{n}$; thus $V_{1}^{+}=\Span(b_{1}, \ldots, b_{r})$, $V_{1}^{-}=\Span(b_{r+1}, \ldots, b_{n})$ and the hermitian form on $V$ has matrix $\ri \cd\twomat {0_{r}} {-1_{r}}{ 1_{r}} {0_{r}}$ in this basis. We may then identify $V_{1}^{\pm}=F^{n}$ and 
 $H=\GL(V_{1})=\GL(W_{1}^{\vee})=\GL_{n}(F)$.

We  will prove the follwoing more precise form of Proposition \ref{int}.
\begin{prop} \lb{inttt2} For $x\in V^{r} $, set $K_{x}^{(v)} := K_{x}\cap \nu^{-1}(C(\vpi_{v}))\subset K_{x}.$
\begin{enumerate}
\item
There  exists a decomposition  into $\Z$-valued Schwartz functions with disjoint supports 
$$\phi^{\bullet} =\phi^{\bullet , \obar} + \phi^{\bullet, \ts}$$ such that:
\begin{itemize}
\item $\phi^{\bullet, \obar}$ takes values in $(q-1)\Z$;
\item if $x\in  {\r{Spt}}(\phi^{\bullet, \ts})$,  then $K_{x}^{(v)}:= K \cap H$ has Galois-level at least~$1$. 
\end{itemize}
\item For all $x\in \r{Spt}(\phi^{\bullet})$, we have $\phi^{\bullet}(x)  [K_{x}:K_{x}^{(v)}]^{-1}\in \Z.$
\end{enumerate}
\end{prop}

We start by studying the stabilisers of the functions $\phi^{\bullet}_{0, y}$. 

\begin{lemm}  \lb{lemma ky} 
 Let $y\in Y\cong  \vpi^{-1}\rM_{r}(\sO_{F})/\rM_{r}(\sO_{F}) $, and let 
 $$K_{y}':= \{d\in \GL_{r}(\sO_{F}) \ | \  (d^{\rt}-1_{r})y\subset \rM_{r}(\sO_{F}) \}.$$
 The stabiliser $K_{y}$ of 
$ \phi_{0, y}^{\bullet}$ under the action of $K^{\circ}\subset H$ is  the subgroup
\beq \lb{ky}
\twomat
 {1_{r} + \vpi\rM_{r}(\sO_{F})}
{\rM_{r}(\sO_{F})}
{\vpi \rM_{r}(\sO_{F})}
{K_{y}'}
  \supset I_{1}.
 \eeq
\end{lemm}
\begin{proof}
It is clear that  $I_{1}\subset K_{y}$, so computing $K_{y}$ is equivalent to computing its image $\ol{K}_{y}\subset \H({k})$.   

Let $h=\twomat abcd \in K^{\circ}$, where  $a, b,c,d\in \rM_{r}(\sO_{F})$. We have
\beq \lb{h act}
h
\twomat  {x_{1, +}}  {x_{2, +}}  {x_{1, -} } {x_{2, -}}
= 
\left(
\twomat abcd {x_{1, +}\choose x_{1, -}}
,
\twomat {a'}{b'}{c'}{d'} {x_{2, +}\choose x_{2, -}}
\right)
\eeq
where the unitarity of $h$ means that $\twomat {a'}{b'}{c'}{d'} $ is characterised by
\beq\lb{h un}
\begin{cases}
a^{\rt} d' - c^{\rt}b'=1_{r}  \\
 d^{\rt}a' - b^{\rt} c' =1_{r} 
 \end{cases} \qquad 
\begin{cases}
   b^{\rt} d'-d^{\rt}b'=0_{r}\\
   c^{\rt}a' + a^{\rt} c' =0_{r}.
 \end{cases}
\eeq
Denote by `$\equiv$' the relation of congruence modulo $\vpi$ on free $\sO_{F}$-modules. 
If $h\in K_{y}$, then for $\displaystyle {x_{1, +}\choose x_{1, -}} ={ 1_{r}\choose 0_{r}} $ we have
  $\displaystyle \twomat abcd {x_{1, +}\choose x_{1, -}}\equiv { 1_{r}\choose 0_{r}} $, 
so that
  $$a\equiv 1_{r}, \qquad c\equiv 0_{r},$$ 
 and \eqref{h un} implies $ \twomat {a'}{b'}{c'}{d'} \equiv  \twomat {d^{\rt, -1}} {d^{\rt, -1} b^{\rt}} {0_{r}} {1_{r}}$. 
  From this,  \eqref{h act}, and the definition of $\phi_{0, y}^{\bullet}$,  we see that $h\in K_{y}$ if and only if further 
$$d^{\rt, -1} y + d^{\rt, -1} b^{t} \rM_{r}(\sO_{F})  \subset y + \rM_{r}(\sO_{F}), $$
that is $(d^{\rt}-1_{r}) y  \subset  \rM_{r}(\sO_{F}) $, as desired. 
\end{proof}

For $y\in Y$, write $d(y)=d$ if the image of $\vpi y$ in $\GL_{r}({k})$ has rank~$d$, and let  
$$Y^{\ts}:=\{y\in Y\ | \ d(y)= r\} \subset Y.$$ 
From Lemma \ref{lemma ky}, we deduce the following.

\begin{lemm}\lb{mult} Suppose $y\in Y-Y^{\ts}$, and let  $K_{y}:=\eqref{ky}$.
The integer $|K_{y}/ I_{1}|$  is a multiple of $q-1$.
\end{lemm}
\begin{proof}  Note that if $d(y)=d$, then the reduction $\ol{K}_{y}'\subset \GL_{r}({k})$ of $K_{y}'$ is $\GL_{r}({k})$-conjugate to $\ol{K}'_{d}:= \twomat {1_{d}} {\rM_{d, r-d}({k})} {0_{r-d, r}}
{\GL_{r-d}({k})}$. 
Therefore $|{K}_{y}/I_{1}| = |\ol{K}_{d}'|$, and when $d<r$ 
 the determinant maps the last group onto $ {k}^{\ts}$. 
  \end{proof}
  
  \begin{rema} \lb{rem d} 
 For $x\in {\mathrm{Spt}} (\phi^{\bullet})$, it is easy to see that there is exactly one integer $0\leq d\leq r$ such that for some $y\in Y$ with $d(y)=d$ and some $h\in K^{\circ}$, we have $h\in {\mathrm{Spt}}(h\phi_{0,y}^{\bullet})$. We denote this integer by $d(x)$. 
\end{rema}

Let us  identify $V\ot_{E} W^{+}=V^{r}$ and denote a typical element by $x$ (rather than $x^{+})$. 
For  $x\in  V^{r}$,  denote $V(x)=\Span(x)^{\perp}\subset V$ (an $r$-dimensional hermitian subspace), and let $H(x)=U(V(x))\subset U(V)=H$. 

We complement the result of  Lemma \ref{mult} by   a lower bound for the Galois-level (Definition \ref{gal level}) 
 the support of $\phi_{0,y}^{\bullet}$ for $y\in Y^{\ts}$.

\begin{lemm} \lb{gal lev}
 Let $y\in Y^{\ts}$
and let $x \in \mathrm{Spt}(\phi_{0, y}^{\bullet})\subset V^{r}$.
The group $K_{x} := K^{\circ}\cap H(x)$ has Galois-level at least $1$.
\end{lemm}
\begin{proof}
The same calculations as in the proof of Lemma \ref{lemma ky}  (with the same notation `$\equiv$')
show that if $h\in K_{x}\subset H(x)$, so that $h x=x$ (in particular $h x\smalltwomat {1_{r}}{}{}{\vpi 1_{r}} \equiv x\smalltwomat {1_{r}}{}{}{\vpi 1_{r}}$), then $h\in K_{y} =I_{1}$.
  Thus $\det h\equiv 1$. 
 \end{proof}

\begin{proof}[Proof of Proposition \ref{inttt2}]
Part 2 follows from part 1: if $x\in \r{Spt}(\phi^{\bullet, \obar})$, it suffices to observe that   $[K_{x}: K_{x}^{(v)}]$ is obviously a divisor of $q-1= [C(1):C(\vpi)]$; whereas if   $x\in \r{Spt}(\phi^{\bullet, \ts})$, we have $[K_{x}:K_{x}^{(v)}]=1$ and $\phi(x)\in \sO$. 

We now prove part 1. 
We have 
\beqq
\phi^{\bullet} &=\sum_{h\in K^{\circ}/I_{1}}\sum_{y\in Y} h\phi^{\bullet}_{0, y}=
\sum_{y\in Y}\sum_{h\in K^{\circ}/K_{y}} |K_{y}/I_{1}|\cd h\phi^{\bullet}_{0, y} \\
 &= \phi^{\bullet, \obar} + \phi^{\bullet, \ts}
 \eeqq
 with 
 \beqq
 \phi^{\bullet, \obar} &:= \sum_{y\in Y-Y^{\ts}} \sum_{h\in K^{\circ}/K_{y}} {|K_{y}/I_{1}|} \cd h\phi^{\bullet}_{0, y} .\\
 \phi^{\bullet, \ts}&:=\sum_{y\in Y^{\ts}} \sum_{h\in K^{\circ}/I_{1}} h\phi_{0, y}^{\bullet}.
\eeqq
By Remark \ref{rem d}, the supports of $ \phi^{\bullet, \obar}$ and  $\phi^{\bullet, \ts}$ are disjoint.  Both functions take values in $\Q$ by the same argument as in Lemma \ref{rat}.2, and both are clearly integral. 

By Lemma \ref{mult}, the function $ \phi^{\bullet, \obar} $ takes values in $(q-1)\Z$. By Lemma \ref{gal lev}, if $x\in \r{Spt}(\phi^{\bullet, \ts}) = \bigcup_{h\in K^{\circ}, y\in Y^{\ts}} \r{Spt} (h\phi_{0, y}^{\bullet})$, then $K_{x}^{(v)}$ has Galois-level at least~$1$.
\end{proof}

\subsection{Test vectors, norm relations and integrality at places in $\wp$}
 
We suppose now that $v\in \wp$.

\subsubsection{P-ordinary representations}
We recall two notion of ordinariness and show that they correspond under Langlands duality.

\begin{defi}  \lb{ord P gal} A de Rham representation $\rho\colon G_{F}\to \GL_{n}(\Qpb)$ is said to be \emph{Panchishkin-ordinary}  (after  \cite[\S~6.7]{Nek93})  if there exists a 
 short exact sequence
$$
0\to \rho^{+}\to \rho\to\rho^{-}\to 0
$$
 of de Rham representations of $G_{F}$ with coefficients in $\Qpb$, such that 
$$
\rF^{0}{\bf D}_{\mathrm{dR}}(\rho^{+})= {\bf D}_{\mathrm{dR}}(\rho^{-})/\rF^{0} {\bf D}_{\mathrm{dR}}(V^{-})=0.
$$
\end{defi}

We denote by  ${\mathrm{Fr}}\in G_{F}$  a lift of the geometric  Frobenius corresponding to the chosen uniformiser $\vpi$. We denote by $\mathrm{WD}(\rho)$  the Weil--Deligne representation over $\Qpb$ attached to a de Rham representation $\rho$  by \cite{Fon94} (see also \cite[\S~1]{TY07}).

\begin{rema}\lb{slopes}
Suppose that $\rho$ is Panchishkin-ordinary, and  put ${\rr}^{\pm}:= \mathrm{WD}(\rho^{\pm})^{\text{Fr-ss}}$, where the superscript denotes Frobenius-semisemplification. By construction, the multiset of slopes ( = $p$-adic valuations of  eigenvalues) of ${\mathrm{Fr}}$  on  $\rr^{\pm}$ coincides with the multiset of slopes of the $[F_{0}:\Q_{p}]^{\text{th}}$ power of the crystalline Frobenius on ${\bf D}_{\text{pst}}(\rho^{\pm})$. In particular, we observe all eigenvalues of  of ${\mathrm{Fr}}$  on  $\rr^{+}$ (respectively $\rr^{-}$) have strictly negative  (respectively non-negative) $p$-adic valuation, and this condition uniquely determines $\rr^{\pm}$ and $\rho^{\pm}$ (up to isomorphism).
\end{rema}

\begin{rema} Suppose that $\rho$ is Panchishkin-ordinary, and assume that  the $\jmath$-Hodge--Tate weight     of $\det \rho^{\pm}$ is independent of $\jmath\colon F\into \bC_{p}$ and equal to $\rw^{\pm}$. 
Let
 $\chi_{\mathrm{cyc}}\colon G_{F}\to \Qpb^{\ts}$ be the cyclotomic character, and let 
\beq\lb{alph}
\alpha:= \chi_{\mathrm{cyc}}^{\rw^{+}}\cd \det\rho^{+}(\mathrm{Fr}).\eeq
Since the  Newton and Hodge polygons of $\mathbf{D}_{\mathrm{pst}}(\rho^{+})$ have the same endpoints,  we have that $\alpha\in \Zpb^{\ts}$.
\end{rema}

 Denote by ${\mathrm{Ind}}_{P}^{G}$ the unitarily  normalised induction, and denote by $\xi_{\pi^{?}}$ the central character of a representation $\pi^{?}$ of a general linear griup. The following definition is adapted from 
\cite{Hid98}.

\begin{defi}\lb{ord P} Let $\pi$ be a 
smooth irreducible  generic representation of $G=\GL_{n}(F)$ with coefficients in $\Qpb$.  Let $\rw^{+}\in\Z_{\leq -1}$, $\rw^{-}\in \Z_{\geq 0}$. We say that $\pi$ is \emph{$P$-ordinary for the Hodge--Tate weights $(\rw^{+}, \rw^{-})$}  if there exists a   $G$-equivariant surjection
\beq\lb{rp}
\rp_{\pi}\colon \mathrm{Ind}_{P}^{G} (\pi^{-}\boxtimes \pi^{+}) \to \pi
\eeq
   for some irreducible admissible representations $\pi^{\pm}$ of $G_{r}:=\GL_{r}(F)$ such that
\beq\lb{xi ord}\xi_{\pi^{+}}|\cd|^{r/2+w^{+}}(\vpi)\in \Zpb^{\ts}.\eeq
\end{defi}

\begin{lemm}\lb{ord ord} Let $\rho\colon G_{F}\to \GL_{n}(\Qpb)$ be a  de Rham representation.
  Let $\pi$ be the smooth irreducible  representation of $\GL_{n}(F)$ over $\Qpb$ that corresponds to $\mathrm{WD}(\rho)$ under the local Langlands correspondence, (re)normalised so that $L(\rho,s)=L(s+1/2, \pi)$. 

Suppose that $\rho$ is Panchishkin-ordinary and the $\jmath$-Hodge--Tate weight  of $\det \rho^{\pm}$ is independent of $\jmath\colon F\into \bC_{p}$ and equal to $\rw^{\pm}$.

Then $\pi$ is $P$-ordinary for the Hodge--Tate weights $(\rw^{+}, \rw^{-})$,  with $\pi^{\pm}$ in \eqref{rp} the representation corresponding to $\mathrm{WD}(\rho^{\pm})$ under the above local Langlands correspondence, and
$$\alpha=|\vpi|^{r/2+\rw^{+}}\xi_{\pi^{+}}(\vpi).$$
\end{lemm}
\begin{proof} 
We freely use the theory of Bernstein--Zelevinsky and the properties of the local Langlands correspondence as summarised for instance in \cite[\S~2]{LLC}, with the notation used there.   We denote by
$$\pi_{-1/2}\colon \rr'\mapsto \pi_{\ru}(\rr(-{1/ 2})), \qquad \pi_{-1/2, \text{ss}}\colon  \rr'\mapsto \pi_{\ru, \text{ss}}(\rr'_{|W_{F}}(-{1/ 2}))$$
the twists of the  unitarily normalised local Langlands correspondence  and, respectively, semisimple local Langlands correspondence of \emph{loc. cit.}

Let $\pi^{\pm}=\pi_{-1/2}(\rr^{\pm})$, and 
write $\pi^{?}=\pi_{\ru}({\bf s}^{?})$ for some multisegments ${\bf s}^{?}$. Since $\pi$ is generic, the segments in ${\bf s}$ can be ordered so as to satisfy the `does not precede' condition above \cite[(2.2.1)]{LLC}. By construction,
  the same  is true of ${\bf s}^{\pm}$, which implies that $\pi^{\pm}$ is generic;  and, together with Remark \ref{slopes},  no  segment in ${\bf s}^{+}$ precedes a segment in ${\bf s}^{-}$. Now (ii) implies that  the unique irreducible quotient of   $\mathrm{Ind}_{P}^{G} (\pi^{-}\boxtimes \pi^{+})$ is generic;  and by construction, its supercuspidal support is the representation $\pi_{-1/2, {\mathrm{ss}}}(\rr)$. But there is a unique up to isomorphism generic irreducible representation with a given supercuspidal support. Since $\pi$ is also  generic and irreducible and has  supercuspidal support  $\pi_{-1/2, {\mathrm{ss}}}(\rr)$, we conclude that the surjection \eqref{rp} exists. The formula for $\alpha$ is clear.
\end{proof}

\subsubsection{The elements $\vphi^{\ra}$ and $f^{\ra}$} 
We specialise back  to our running assumptions, so that $\rho$ is pure of weight~$-1$ and Panchishkin-ordinary with $\jmath$-Hodge--Tate weights $\{-r, \ldots, r-1\}$ (for every $\jmath\colon F\into \bC_{p}$); and $\pi$ is the associated $P$-ordinary representation of $G$. In the notation of the previous paragraphs, we have $\rw^{+}=-{r\choose 2}$, $\rw^{-}= {r-1 \choose 2}$. 
Then  $\pi^{\vee}$ is also $P$-ordinary for the weight $(\rw^{+}, \rw^{-})$ with respect to the representations $\pi^{\vee, \pm}:=\pi^{\mp, \vee}$.
 We note that $\alpha=q^{r^{2}/2}\xi_{\pi^{+}}(\vpi)$ and put
$$ \alpha^{\vee}:= q^{r^{2}/2}\xi_{\pi^{\vee, +}}(\vpi)= \xi_{\pi}^{-1}(\vpi) \alpha \in \Zpb^{\ts}.$$
 We fix perfect parings $(\ , \ )_{\pi^{\pm}}\colon\pi^{\vee, \mp}\ot\pi^{\pm}\to \Qpb$.

For any $s\in \Z$, we define elements of $G$ by
$$t:= \twomat{\vpi 1_{r}}{}{}{1_{r}}, \qquad w_{s}:=\twomat  {}{1_{r}} {-\vpi^{s}1_{r}}{} =w t^{s}$$
and we put $t^{w}=w^{-1}tw=\smalltwomat {1_{r}}{}{}{\vpi 1_{r}}$. We also put
$$U_{t}=\sum_{b\in \rM_{r}(k)} \twomat {1_{r}} {b}{}{1_{r}} t \qquad \in \Z[G];$$
then 
for all $s'\geq s\geq 1$, the double coset operator
$$I_{s'}tI_{s}$$
acts by $U_{t}$ on any smooth $G$-module.

We introduce the following condition, which we assume from now on:
\beq\lb{cond: marcil}
\text{either $\pi$ is unramified or $\pi^{\pm}$  are both supercuspidal.}
\eeq
\begin{rema}\lb{rema: marcil} The reason for imposing this  technical condition is that it(s second part) appears in the current version of \cite{Mar1}; it is expected to be removed in  future updates.
\end{rema}

\begin{prop} \lb{mat coeff p}
There exist 
$$\vphi^{\pm}\in \pi^{\vee, \pm}, \quad f^{\pm}\in  \pi^{\pm}, \quad \vphi^{\ord}\in \pi, \quad f^{\ord}\in \pi$$
  satisfying the following conditions: 
\begin{enumerate}
\item  $(\vphi^{-}, f^{+})_{\pi^{+}}=(\vphi^{+}, f^{-})_{\pi^{-}}=1$;
\item
\lb{Up eigen}
 there is a constant  $c(\pi)$ such that the vectors $\vphi^{\ord}$, $f^{\ord}$  are invariant under $I_{c(\pi)}$  and satisfy
\beqq
U_{t} \vphi^{\ord} =\alpha^{\vee} \vphi^{\ord}, \qquad
U_{t}   f^{\ord} = \alpha f^{\ord};
\eeqq
\item \lb{inn prod D}
 for each $c\geq c(\pi)$, setting 
\beqq
\vphi_{c}^{\ra}:= q^{cr^{2}}\alpha^{\vee,-c} \pi^{\vee}(w_{c})\vphi^{\ord}, 
\qquad 
f_{c}^{\ra}:= q^{cr^{2}}\alpha^{-c} \pi(w_{c})f^{\ord},\eeqq
we have, for every $y_{\pm}\in \GL_{r}(F)$, 
\beqq
( \pi^{\vee}(\smalltwomat {y_{+}} {}{}{y_{-}} w^{-1}) \vphi_{c}^{\ra}, f_{c}^{\ra})_{\pi}
= 
| \det y_{+}|^{r/2} | \det y_{-}|^{-r/2}\cd (\pi^{\vee, +}(y_{+}) \vphi^{+}, f^{-})_{\pi^{-}}  ( \pi^{\vee, -}(y_{-})\vphi^{-},  f^{+})_{\pi^+}.
\eeqq
\end{enumerate}
\end{prop}

(The superscripts `a' stand for `anti-ordinary'.)

\begin{proof}

Consider first the case where $\pi^{\pm}$ are supercuspidal.
We deduce the proposition  form the  results of  Marcil in \cite{Mar1, Mar2}. 
Let $c(\pi)\geq 1$ be  the  minimal integer (denoted by $r$ in \cite{Mar1, Mar2}) for which the construction that we are about to  cite can be performed.  
  Let
$$f^{+},\ f^{-}; \quad \vphi^{-},\ \vphi^{+}; \quad f_{c(\pi)}^{\rM} \in \pi^{I_{c(\pi)}},\ \vphi_{c(\pi)}^{\rM}\in (\pi^{\vee})^{I_{c(\pi)}^{\rt}}$$
 be the vectors denoted respectively 
$$\phi_{a}, \ \phi_{b}; \quad \mu \wt\phi_{a},  \mu \wt\phi_{b};\quad \vphi, \ \wt{\phi}$$
 in \cite[pp. 14-15]{Mar2}  (we omit all the subscripts `$w$'  when transcribing notation from \cite{Mar1, Mar2}), where we adjust the scalar $\mu\in \Q^{\ts}$ so that $(\vphi^{\mp}, f^{\pm})_{\pi^{\pm}}=1$. 
As noted in \cite[Remarks 2.8, 2.11]{Mar2}, the existence  of such vectors, which may depend on some choices that we fix, is guaranteed under our assumptions by \cite[Theorem 4.3, Lemma 4.6]{Mar1}. 

Moreover, denote by $\delta_{P}\colon\smalltwomat {g_{+}} {}{}{g_{-}}\mapsto |\det g_{+}|^{r}|\det g_{-}|^{-r}$ the modulus character of $P$; then  by those results we have 
\beqq \lb{Up e1}
I_{c(\pi)}t^{-1}I_{c(\pi)}f_{c(\pi)}^{\rM} = \delta_{P}^{-1/2}(t) \xi_{\pi^{\vee, +}}(\vpi) f_{c(\pi)}^{\rM},\eeqq
 that is
\beq \lb{Up e1}
I_{c(\pi)} t^{w}I_{c(\pi)}f_{c(\pi)}^{\rM}= \delta_{P}^{-1/2}(t) \xi_{\pi^{\vee, -, \vee}}(\vpi) f_{c(\pi)}^{\rM}
= \alpha f_{c(\pi)}^{\rM};  \eeq
 and similarly
\beq \lb{Up e2}
I^{\rt}_{c(\pi)}t I^{\rt}_{c(\pi)}\vphi_{c(\pi)}^{\rM}=   \alpha^{\vee}\vphi_{c(\pi)}^{\rM}.\eeq
 
 Finally, by the proof of \cite[Proposition 4.3]{Mar2}:
\beq\lb{inn prod M}
( \pi^{\vee}(\smalltwomat {g_{+}} {}{}{g_{-}}) \vphi_{c}^{\rM}, f_{c}^{\rM})_{\pi}
= 
\mu_{c}'\cd
| \det g_{+}|^{r/2} | \det g_{-}|^{-r/2}\cd (\pi^{\vee, +}(g_{+}) \vphi^{+}, f^{-})_{\pi^{-}}  ( \pi^{\vee, -}(g_{-})\vphi^{-},  f^{+})_{\pi^+},
\eeq
for some constant $\mu'_{c}$.

Let
\beqq f_{c(\pi)}^{\ra}:= \mu'_{c(\pi)}{}^{-1}\cd f_{c(\pi)}^{\rM}, &\qquad \vphi_{c(\pi)}^{\ra}:= \pi^{\vee}(w) \vphi_{c(\pi)}^{\rM}, \\
f^{\ord}:=w_{c(\pi)}^{-1}f_{c(\pi)}^{\ra}, & \qquad\vphi^{\ord}:=w_{c(\pi)}^{-1}\vphi_{c(\pi)}^{\ra};
\eeqq
then by the definitions and \eqref{Up e1}, \eqref{Up e2}
we have part \ref{Up eigen}; and by \eqref{inn prod M} we have
  part  \ref{inn prod D}.

The case where $\pi $ is unramified can be  similarly deduced from \cite[Proposition 3.25]{DL}, by 
 taking $c(\pi)=1$ and  $\vphi_{1}^{\ra}=w_{1} \vphi^{\ord}$,  $f_{1}^{\ra}=w_{1}f^{\ord}$   to be the vectors denoted respectively  by $\vphi_{v}^{\vee}$ and $\vphi_{v}$ in \emph{loc. cit.}
  \end{proof}

\subsubsection{Schwartz functions} Let us first introduce some notation.
For a compact open subgroup  $J\subset \GL_{r}(F)$ and any Haar measure $\rd y$ on $J$, put
\beqq
\delta_{J}(y)&:= \vol(J, \rd^{}y)^{-1}\cd \one_{J}(y)\\
\psi_{J}(y')&:=\vol(J ,\rd^{}y)^{-1}\cd \int_{J}\psi(\Tr(y^{\rt}y')) \,\rd^{}y.\eeqq
If $\Phi=(\Phi_{ij}) $ is a matrix of  $\Qpb$-valued functions of the variables $X_{ij}$ (viewed as the entries of a matrix $X$),
we write 
$$\Phi(X):= \prod_{i,j}\Phi_{ij}(X_{ij}).$$

Let $J_{+}\subset \GL_{r}(\sO_{F})$, respectively  $J_{-}\subset \GL_{r}(\sO_{F})$, be an open  subgroup fixing $\vphi^{+}$ and $f^{-}$, respectively $\vphi^{-}$ and $f^{+}$. We fix the generator $d=\vpi^{v(d)}\in\sO_{F}$ of the different ideal to be a power of $\vpi$.
For every pair of integers $s,s'\geq 0$, let\footnote{For convenience,  the choices in our  naming of the coordinates here are different from  those of \S~\ref{sec: horiz}.}
 \beq\lb{phis}
  \phi^{(s, s')}
 \left(\twomat  {x }  {x'}  {y} {y' }\right)
& :=
\one[
 x^{\rt}y' - y^{\rt}x'\in d^{-1}\rM_{r}(\sO)]
\cd
 \twomat
  {\one_{\vpi^{-s}d^{-1}\rM_{r}(\sO)}}
 {  \one_{\vpi^{-s'}\rM_{r}(\sO)}}
{\delta_{\vpi^{-s}d^{-1}J_{-}}}
{\delta_{\vpi^{-s'}J_{+}}}
\left(\twomat  {x}  {x'}  {y } {y'}\right)
\eeq
\begin{rema}
\lb{comp DL schw} This  is a variant of the function $\phi_{2,v}^{[(s,s')_{v}]}$ defined before \cite[Lemma 4.29]{DL}; more precisely, let $m_{d}:=\twomat {d 1_{r}} {}{}{1_{r}}$; then, when 
 $J_{\pm}=\GL_{r}(\sO_{F})$, we have $\phi^{(s,s')}=m_{d}\phi_{v,2}^{[(s,s')_{v}]}$.
\end{rema}

For the next two lemmas, recall we have fixed an isomorphism between the isomorphic groups $G$ and $H$, but that they act rather differently on $\sS_{\chi}(V^{r})$; thus we add a superscript `$G$' in the notation for $U_{t}\in \Z[G]$, and denote by $U_{t}^{H}\in \Z[H]$ the corresponding operator for $H$; likewise for the subgroups $I_{c}^{?}$ of $G$ and $H$. 
\begin{lemm} \lb{Up phi}
 There exists an integer $c\geq 1$ such that for all $s,s'\geq0$, $\phi^{(s,s')}$ is fixed under $I_{c}^{G}\ts I_{c}^{H}$, and 
\beqq
\omega_{\chi}(U_{t}^{G}) \phi^{(s, s')}&=\chi(\vpi)^{-r}  \phi^{(s+1, s')},\\
\omega_{\chi}( U_{t}^{H}) \phi^{(s, s')}&=\chi(\vpi)^{r} \phi^{(s, s'+1)}. \eeqq
\end{lemm}
\begin{proof} 
It suffices to check that $\phi^{(0,0)}$ is fixed under $\rN(\sO_{F})^{G}\ts \rN(\sO_{F})^{H}\subset G\ts H$, which is straightforward, and the two formulas. 

The first formula is verified   as in  \cite[Lemma 4.29 (1)]{DL}; to compare, see Remark \ref{comp DL schw}, and note that our $U_{t}^{G}$ equals the operator $m_{d}\rU_{v}^{1_{w_{1}}}m_{d}^{-1}$ used in   \emph{loc. cit.}
 
 For the second formula,  by the definitions and \eqref{h un}, we have 
$$U_{t}^{H}\phi\left({x \choose y}, {x'\choose y'}\right)
=\chi(\vpi)^{r} \sum_{b\in \rM_{r}(\sO_{F}/\vpi)} \phi\left(\twomat {\vpi^{-1}1_{r}} {\vpi^{-1}b}{}{1_{r}} {x\choose y},
\twomat {1_{r}} {b^{\rt}} {}{\vpi 1_{r}} {x'\choose y'} \right).
$$
For   $\phi=\phi^{(s,s')}$, it is easy to check that each term in the sum vanishes unless $\twomat x{x'} y {y'}$ belongs to the support of $\phi^{(s, s'+1)}$, whereas if this condition is satisfied, then the sum contains only one nonzero term, with value $\vol(J_{+})^{-1}\vol(J_{-})^{-1}$, which is the one indexed by the class 
$$b=- xy^{-1} \equiv (x'y'^{-1})^{\rt} \in \rM_{r}(\sO_{F}/\vpi).$$
\end{proof}

\begin{defi}\lb{test p}
Let $c(\pi)$ and $\vphi^{\ra}_{c}$, $f^{\ra}_{c}$ be as in  Proposition \ref{mat coeff p}, and let  $c\geq c(\pi)$ be the minimal integer satisfying the conditions of Lemma \ref{Up phi}. We define for all $s\geq 0$:
\beqq
\vphi^{\ra}:= \vphi^{\ra}_{c} \quad \in \pi, 
\qquad 
\phi^{(s)}&:=\alpha^{-s-v(d)}\alpha^{\vee, -s}\cd \phi^{(s,s)} \quad \in \sS(V^{r}),
\qquad
f^{\ra}:= f^{\ra}_{c} \quad \in \sg,\\
 \lm^{(s)}&:=\vphi^{\ra}\ot \phi^{(s)}\ot f^{\ra}.
 \eeqq
\end{defi}

\subsubsection{Norm relations}
We prove the local form of the vertical norm relations.

\begin{lemm}\lb{p pet}
For every smooth admissible $G\ts H$-module $\sS$, every $1\leq d\leq c$, 
and every $\phi\in\sS^{I_{c}^{G}\ts I_{c}^{H}}$, we have
\beqq
[\vphi^{\ra}\ot U_{t}^{G}
\phi]&= \alpha^{}\cd [\vphi^{\ra}\ot \phi]  & \text{in } (\pi^{\vee}\ot\sS)_{G},\\
[U_{t}^{H}\phi \ot f^{\ra}]&= \alpha^{\vee}\cd [\phi \ot f^{\ra}]  & \text{in } (\sS\ot \sg)_{H}.
\eeqq
\end{lemm}
\begin{proof}
For the first equality, dropping all superscripts $G$, we have
\beqq
[\vphi^{\ra}\ot I_{d}tI_{d} \phi] 
&=
[\pi^{\vee}(I_{d}t^{-1} I_{d} w_{c}) \vphi\ot \phi^{\ord}]
=
[\pi^{\vee}(I_{c}t^{-1} I_{c} w_{c})\vphi^{\ord}\ot \phi]\\
&=
[\pi^{\vee}(w_{c}I_{c}t^{w,-1}I_{c} )\vphi^{\ord}\ot \phi]
=  \xi_{\pi^{\vee}}(\vpi)^{-1}
\alpha^{\vee}\cd
[ \pi^{\vee}(w_{c})\vphi^{\ord}\ot \phi]
=\alpha\cd
[ \vphi^{\ra}\ot \phi].
\eeqq
The proof of the second equality is virtually identical.
\end{proof}
\begin{prop}[\emph{\ $=$\, Proposition \ref{prop NRv}}] For every open compact subgroup $C\subset F^{\ts}$, the image of $\lm^{(s)}$ in $\Lm_{\rho}^{C}$ is independent of $s\geq 0$. 
\end{prop}
\begin{proof} This follows from Definition \ref{test p}, Lemma \ref{Up phi}, and Lemma \ref{p pet} applied to $\sS=\sS_{\chi}(V^{r})$ for each character $\chi$ of $F^{\ts}/C$. 
\end{proof}

\Subsubsection{Integrality}
\begin{prop} \lb{intt p}
For every $s\geq 0$, we have 
$$\phi^{(s)}\in \sS(V^{r}, \ol{\mathbf{Z}}_{p, C(\vpi^{s})})^{I_{c}\ts I_{c}}.$$
\end{prop}
\begin{proof}
It suffices to prove that for each $x\in \mathrm{Spt}(\phi^{(s)})$, the group $K_{x}$ has Galois-level at least $s$ (Definition \ref{gal level}). 
  This is proved in \cite[Lemma 4.36]{DL} (for a slightly different Schwartz function, but the same proof goes through).
\end{proof}

\begin{proof}[Proof of Proposition \ref{int}] It follows from Propositions \ref{inttt2}.2 and \ref{intt p}.
\end{proof}

\subsubsection{Local non-vanishing} \lb{loc nv}
We conclude by studying the image of $\lm_{}^{(s)}$ in $\Lm_{\rho}$. 
\begin{prop} \lb{gamma}
For every $s\geq 0$, we have
 $$\tht(\lm^{(s)})
 = q^{v(d)r^{2}/2} \gamma(1/2, \pi^{\vee, +}, \psi)^{-1}.
 $$
\end{prop}
\begin{proof}[Proof of Proposition \ref{NV}, assuming Proposition \ref{gamma}]
We restore the notation used in the global context. Denote by $[\lm_{v}^{(0)}] $ the image of $\lm_{v}^{(0)}$  in $\Lm_{\rho,v}$. It is clear that $[\lm_{v}^{(0)}]\neq 0$ at all $v\in S$. At $v\notin S\wp$, the nontriviality of $[\lm_{v}^{(0)}]$ follows from \cite[Proposition 3.6.4]{cet}.  At $v\in \wp$,  the desired assertion follows from  Proposition \ref{gamma}  since $\gamma(1/2+s, \pi_{v}^{\vee, +}, \psi_{w})
=\gamma(\mathrm{WD}(\rho_{{w_{2}}}^{+}), \psi_{w},s).$
 \end{proof}

In order to prove Proposition \ref{gamma}, we need a lemma. Let $\phi':=\sF\phi^{(0)}$, and consider the map  
\beq\lb{gabcd}
&g\colon \GL_{r}(F)\ts \rM_{r}(F)\ts \rM_{r}(F)\ts \GL_{r}(F)\longrightarrow G=\GL_{n}(F)\\
 &g(y_{+},x_{+},x_{-},y_{-}):=\twomat 1{x_{+}}{}1 \twomat {y_{+}}{}{}{y_{-}} \twomat 1{}{x_{-}}1 w^{-1}.
\eeq
 \begin{lemm} We have 
 $$\phi'(g(y_{+},x_{+},x_{-},y_{-}))
= |d|^{r^{2}/2}
\cd
 \twomat 
{\psi_{J_{+}}}
{\one_{\rM_{r}(\sO)}}
{\one_{\rM_{r}(\sO)}}
{\delta_{d^{-1} J_{-}}}
\left(
\twomat
{y_{+}}
{x_{+}}
{x_{-}}
{y_{-}}
\right).
$$
\end{lemm}
\begin{proof}
With the change of variables
$x_{ +}' = y_{ -}^{\rt, -1}x_{ +}^{\rt}y'_{-}+x_{+}''$, we  have 
\beqq
\phi'\left(\twomat {x_{+}} {y_{+}} {y_{-}} {x_{-}} \right)&=
\int_{\rM_{r}(F)}\int_{\rM_{r}(F)}
\phi\left(\twomat  {x_{+}} {x'_{+}} {y_{-}} {y'_{-}} \right)\, 
\psi(x'_{ +}{}^{\rt}x_{-})\, \psi(-y'_{ -}{}^{\rt}y_{+}) 
\,  \rd x'_{+} \rd y'_{ -}  \\
&= \one_{
d^{-1}\rM_{r}(\sO)}(x_{+})\delta_{
d^{-1}J_{-}} (y_{-}) 
\int_{
\rM_{r}(F)}
\int_{
\rM_{r}(\sO)}
\psi(( y_{-}'{}^{\rt}x_{+}y_{-}^{-1}+  x_{+}'' {}^{\rt} )x_{-}) \, 
\delta_{J_{+}}(y'_{-})
\,  \psi(-y'_{-}{}^{\rt}y_{+} ) \,  \rd x_{+}'' \rd y_{ -}' 
\\
&=  \one_{
d^{-1} \rM_{r}(\sO)}(x_{+})\delta_{
d^{-1} J_{-}} (y_{-})
\int_{
\rM_{r}(\sO)}\psi(  x_{+}''{}^{\rt}x_{-})\,  \rd x_{+}'' 
\int_{\rM_{r}(F)}
\delta_{J_{+}}(y_{-}')  \psi(y_{-}'{}^{\rt}(x_{+}y_{+}^{-1} x_{-} -y_{+})) \rd y_{-}' \\
 &= 
 |d|^{r^{2}/2}
  \one_{
d^{-1}  \rM_{r}(\sO)}(x_{+})\delta_{
d^{-1}   J_{-}}(y_{-})
  \one_{
d^{-1}  \rM_{r}(\sO)}(x_{-})
  \psi_{J_{+}}
(x_{+}y_{-}^{-1} x_{-} -y_{+})\\
 &=  |d|^{r^{2}/2}
 \twomat 
{\one_{
d^{-1} \rM_{r}(\sO)}}
{\psi_{J_{+}}}
{\delta_{d^{-1}  J_{-}}}
{\one_{d^{-1} M_{r}(\sO)}}
\left(
\twomat
{x_{+}}
{x_{+}y_{-}^{-1} x_{-} -y_{+}}
{y_{-}}
{x_{-}}
\right).
  \eeqq
Then the desired formula follows from  evaluating at 
$$g(y_{+},x_{+},x_{-},y_{-})= 
\twomat{x_{+}y_{-}} {-y_{+}-\, x_{+}y_{-}x_{-}} {y_{-}}{-y_{-}x_{-}}.$$
\end{proof}

\begin{proof}[Proof of Proposition \ref{gamma}]
By the definitions, we have
$$
\tht(\lm^{(s)})= \tht(\lm^{(0)})
=\int_{G} (g\vphi_{c}^{\ra}, f_{c}^{\ra})_{\pi}\cd \phi'(g)
     |\det g|^{r}
\, dg.$$
We integrate over the full-measure subset that is the image of the map $g=\eqref{gabcd}$, for which  
 $$\rd g(y_{+},x_{+},x_{-},y_{-})= |\det y_{+}|^{-r}|\det y_{-}|^{r} \rd x_{-} \rd^{}y_{+} \rd^{}y_{-} \rd x_{+}.$$
where the Haar measures $\rd y_{\pm}$ on $G_{r}=\GL_{r}(F)$ and $\rd x_{\pm}$ on $\rM_{r}(F)$
   are normalised by assigning volume $|d|^{r^{2}/2}$ respectively to 
    $\GL_{r}(\sO_{F})$ and $\rM_{r}(\sO_{F})$. 
    Then we obtain 
\beqq
\tht(\lm^{(0)})
&= \int_{G}
\left(\pi^{\vee}\left(
g(y_{+},x_{+},x_{-},y_{-})\right)\vphi^\ra, f^{\ra}\right)_{\pi}  
\cd \phi'
\left(
g(y_{+},x_{+},x_{-},y_{-})
\right)\,
 |\det y_{+} \det y_{-}|^{r} 
\, \rd g(y_{+},x_{+},x_{-},y_{-}) \\
&
\begin{multlined}
=\alpha^{-v(d)} |d|^{r^{2}/2}
\int_{\rM_{r}(\sO_{F})}
\int_{\GL_{r}(F)}
\int_{\rM_{r}(\sO_{F})}
\int_{\GL_{r}(F)}
\left(\pi^{\vee}\left(\smalltwomat 1{x_{+}}{}1 \smalltwomat {y_{+}}{}{}{y_{-}}  w^{-1}\smalltwomat 1{-x_{-}}{}1 \right)\vphi^{\ra}, f^{\ra}\right)_{\pi} 
\\ \cdot 
\psi_{J_{+}}(y_{+})  \delta_{d^{-1}J_{-}}(y_{-})
|\det y_{-}|^{2r}\,  \rd x_{-} \rd^{}y_{+} \rd^{}y_{-} \rd x_{+}.
\end{multlined}
\eeqq
Since $\vphi^{\ra}$ and $f^{\ra}$ are both invariant under $\rN(\sO)$, the integrations in $\rd x_{\pm}$ give $1$, and we get
$$ 
\tht(\lm^{(0)})=\alpha^{-v(d)} |d|^{r^{2}/2}
\int_{\GL_r(F)} \int_{ \GL_r(F)}\left(\pi^{\vee}\left(\smalltwomat {y_{+}}{}{}{y_{-}}  w^{-1}\right) \vphi^{\ra}, f^{\ra}\right)_{\pi} 
\cdot 
\psi_{J_{+}}(y_{+})   \delta_{d^{-1}J_{-}}(y_{-})
|\det y_{-}|^{2r}\, \rd y_{-}  \rd y_{+}  
.$$
By the formula for the matrix coefficient in Proposition \ref{mat coeff p}, we deduce  
$$\tht(\lm^{(0)})=Z_{+} Z_{-},$$ 
where
\beqq
Z_{+}&:=|d|^{-r^{2}/2}\int_{\GL_r(F)}(\pi^{\vee, +}(y_{+})\vphi^{+}, f^{-})_{\pi^{-}} \cdot 
\psi_{J_{+}}(y_{+}) 
  |\det y_{+}|^{r/2}\,  \rd y_{+} , \\
Z_{-}&:=\alpha^{-v(d)} |d|^{r^{2}}
 \int_{\GL_r(F)}
(\pi^{\vee, -}(y_{-})\vphi^{-}, f^{+})_{\pi^{+}} \delta_{d^{-1}J_{-}}(y_{-})
 |\det y_{-}|^{3r/2}\,  \rd y_{-}.
\eeqq

Since $\alpha= |\vpi|^{r^{2}/2} \xi_{\pi^{+}}(\vpi)$, we have 
$$Z_{-}= |d|^{r^{2}} \alpha^{-v(d)}\xi_{\pi^{\vee, -}}(d^{-1})|d|^{3r^{2}/2}=1;$$
 whereas by the Godement--Jacquet  functional equation (\cite[Proposition 1.2 (3)]{Jac79}, which has a typo corrected in (1.3.7) \emph{ibid.}), 
\beqq
|d|^{r^{2}/2} \cd Z_{+}&= \gamma(1/2, \pi^{\vee, +}, \psi
 )^{-1}\, 
\int_{\GL_{r}(F)}   (\vphi^{+}, \pi^{-}(y_{+})f^{-})_{\pi^{-}} \delta_{J_{+}}(y_{+})
   |\det y_{+}|^{r/2} \, \rd y_{+}
=
\gamma(1/2, \pi^{\vee, +} , \psi
)^{-1}.
\eeqq
This completes the proof.
\end{proof}

\end{document}